\documentclass[12pt]{article}
\usepackage[centertags]{amsmath}
\usepackage{amsfonts}
\usepackage{amssymb}
\usepackage{amsthm}
\usepackage{graphics,graphicx}
\usepackage{newlfont}

\newlength{\defbaselineskip}
\newcommand{\setlinespacing}[1]%
           {\setlength{\baselineskip}{#1 \defbaselineskip}}

\newcommand{\cW}{{\cal W}}
\newcommand{\cV}{\mathcal{V}}
\newcommand{\cR}{\mathcal{R}}

\newcommand{\e}{\epsilon}
\newcommand{\de}{\delta}
\newcommand{\al}{\alpha}
\newcommand{\be}{\beta}
\newcommand{\bx}{\bar{x}}
\newcommand{\bt}{\bar{t}}

\newcommand{\cF}{\mathcal {F}}

\theoremstyle{plain}
{\swapnumbers
\newtheorem{thm}{Theorem}[section]
\newtheorem{cor}[thm]{Corollary}

\newtheorem{prop}[thm]{Proposition}
}
\theoremstyle{definition}
{\swapnumbers

\newtheorem{parag}[thm]{}
\newtheorem*{parstr}{\addtocounter{thm}{1}\thesection.\arabic{thm}*}
\newtheorem*{intropar}{\addtocounter{thm}{1}\arabic{thm}}
}
\theoremstyle{remark} {\swapnumbers

}
\numberwithin{equation}{section}

\begin{document}
\title{The entropy formula for the Ricci flow \\ and its geometric
applications }\author{Grisha Perelman\thanks{St.Petersburg branch
of Steklov Mathematical Institute, Fontanka 27, St.Petersburg
191011, Russia. Email: perelman@pdmi.ras.ru or
perelman@math.sunysb.edu ; I was partially supported by personal
savings accumulated during my visits to the Courant Institute in
the Fall of 1992, to the SUNY at Stony Brook in the Spring of
1993, and to the UC at Berkeley as a Miller Fellow in 1993-95. I'd
like to thank everyone who worked to make those opportunities
available to me.}} \maketitle
\section*{Introduction} \begin{intropar}The Ricci flow equation, introduced
by Richard Hamilton [H 1], is the evolution equation
$\frac{d}{dt}g_{ij}(t)=-2R_{ij}$ for a riemannian metric
$g_{ij}(t).$ In his seminal paper, Hamilton proved that this
equation has a unique solution for a short time for an arbitrary
(smooth) metric on a closed manifold. The evolution equation for
the metric tensor implies the evolution equation for the curvature
tensor of the form $Rm_t=\triangle Rm +Q,$ where $Q$ is a certain
quadratic expression of the curvatures. In particular, the scalar
curvature $R$ satisfies $R_t=\triangle R+2|\mbox{Ric}|^2,$ so by
the maximum principle its minimum is non-decreasing along the
flow. By developing a maximum principle for tensors, Hamilton [H
1,H 2] proved that Ricci flow preserves the positivity of the
Ricci tensor in dimension three and of the curvature operator in
all dimensions; moreover, the eigenvalues of the Ricci tensor in
dimension three and of the curvature operator in dimension four
are getting pinched pointwisely as the curvature is getting large.
This observation allowed him to prove the convergence results: the
evolving metrics (on a closed manifold) of positive Ricci
curvature in dimension three, or positive curvature operator in
dimension four converge, modulo scaling, to metrics of constant
positive curvature. \par Without assumptions on curvature the long
time behavior of the metric evolving by Ricci flow may be more
complicated. In particular, as $t$ approaches some finite time
$T,$ the curvatures may become arbitrarily large in some region
while staying bounded in its complement. In such a case, it is
useful to look at the blow up of the solution for $t$ close to $T$
at  a point where curvature is large (the time is scaled with the
same factor as the metric tensor). Hamilton [H 9] proved a
convergence theorem , which implies that a subsequence of such
scalings smoothly converges (modulo diffeomorphisms) to a complete
solution to the Ricci flow whenever the curvatures of the scaled
metrics are uniformly bounded (on some time interval), and their
injectivity radii at the origin are bounded away from zero;
moreover, if the size of the scaled time interval goes to
infinity, then the limit solution is ancient, that is defined on a
time interval of the form $(-\infty , T).$ In general it may be
hard to analyze an arbitrary ancient solution. However, Ivey [I]
and Hamilton [H 4] proved that in dimension three, at the points
where scalar curvature is large, the negative part of the
curvature tensor is small compared to the scalar curvature, and
therefore the blow-up limits have necessarily nonnegative
sectional curvature. On the other hand, Hamilton [H 3] discovered
a remarkable property of solutions with nonnegative curvature
operator in arbitrary dimension, called a differential Harnack
inequality, which allows, in particular, to compare the curvatures
of the solution at different points and different times. These
results lead Hamilton to certain conjectures on the structure of
the blow-up limits in dimension three, see [H 4,$ \S 26$]; the
present work confirms them.
\par The most natural way of forming a singularity in finite time
is by pinching an (almost) round cylindrical neck. In this case it
is natural to make a surgery by cutting open the neck and gluing
small caps to each of the boundaries, and then to continue running
the Ricci flow. The exact procedure was described by Hamilton [H
5] in the case of four-manifolds, satisfying certain curvature
assumptions. He also expressed the hope that a similar procedure
would work in the three dimensional case, without any a priory
assumptions, and that after finite number of surgeries, the Ricci
flow would exist for all time $t\to\infty,$ and be nonsingular, in
the sense that the normalized curvatures
$\tilde{Rm}(x,t)=tRm(x,t)$ would stay bounded. The topology of
such nonsingular solutions was described by Hamilton [H 6] to the
extent sufficient to make sure that no counterexample to the
Thurston geometrization conjecture can occur among them. Thus, the
implementation of Hamilton program would imply the geometrization
conjecture for closed three-manifolds. \par In this paper we carry
out some details of Hamilton program. The more technically
complicated arguments, related to the surgery, will be discussed
elsewhere. We have not been able to confirm Hamilton's hope that
the solution that exists for all time $t\to\infty$ necessarily has
bounded normalized curvature; still we are able to show that the
region where this does not hold is locally collapsed with
curvature bounded below; by our earlier (partly unpublished) work
this is enough for topological conclusions.
\par Our present work has also some applications to the
Hamilton-Tian conjecture concerning $\mbox{K\"{a}hler-Ricci}$ flow
on $\mbox{K\"{a}hler}$ manifolds with positive first Chern class;
these will be discussed in a separate paper.
\end{intropar}
\begin{intropar} The Ricci flow has also been discussed in quantum field
theory, as an approximation to the renormalization group (RG) flow
for the two-dimensional nonlinear $\sigma$-model, see [Gaw,$\S 3$]
and references therein. While my background in quantum physics is
insufficient to discuss this on a technical level, I would like to
speculate on the Wilsonian picture of the RG flow. \par In this
picture, $t$ corresponds to the scale parameter; the larger is
$t,$ the larger is the distance scale and the smaller is the
energy scale; to compute something on a lower energy scale one has
to average the contributions of the degrees of freedom,
corresponding to the higher energy scale. In other words,
decreasing of $t$ should correspond to looking at our Space
through a microscope with higher resolution, where Space is now
described not by some (riemannian or any other) metric, but by an
hierarchy of riemannian metrics, connected by the Ricci flow
equation. Note that we have a paradox here: the regions that
appear to be far from each other at larger distance scale may
become close at smaller distance scale; moreover, if we allow
Ricci flow through singularities, the regions that are in
different connected components at larger distance scale may become
neighboring when viewed through microscope.
\par Anyway, this connection between the Ricci flow and the RG
flow suggests that Ricci flow must be gradient-like; the present
work confirms this expectation.
\end{intropar}
\begin{intropar}
The paper is organized as follows. In $\S 1$ we explain why Ricci
flow can be regarded as a gradient flow. In $\S 2,3$ we prove that
Ricci flow, considered as a dynamical system on the space of
riemannian metrics modulo diffeomorphisms and scaling, has no
nontrivial periodic orbits. The easy (and known) case of metrics
with negative minimum of scalar curvature is treated in $\S 2;$
the other case is dealt with in $\S 3,$ using our main
monotonicity formula (3.4) and the Gaussian logarithmic Sobolev
inequality, due to L.Gross. In $\S 4$ we apply our monotonicity
formula to prove that for a smooth  solution on a finite time
interval, the injectivity radius at each point is controlled by
the curvatures at nearby points. This result removes the major
stumbling block in Hamilton's approach to geometrization. In $\S
5$ we give an interpretation of our monotonicity formula in terms
of the entropy for certain canonical ensemble. In $\S 6$ we try to
interpret the formal expressions , arising in the study of the
Ricci flow, as the natural geometric quantities for a certain
Riemannian manifold of potentially infinite dimension. The
Bishop-Gromov relative volume comparison theorem for this
particular manifold can in turn be interpreted as another
monotonicity formula for the Ricci flow. This formula is
rigorously proved  in $\S 7;$ it may be more useful than the first
one in local considerations. In $\S 8$ it is applied to obtain the
injectivity radius control under somewhat different assumptions
than in $\S 4.$ In $\S 9$ we consider one more way to localize the
original monotonicity formula, this time using the differential
Harnack inequality for the solutions of the conjugate heat
equation, in the spirit of Li-Yau and Hamilton. The technique of
$\S 9$ and the logarithmic Sobolev inequality are then used in $\S
10$ to show that Ricci flow can not quickly turn an almost
euclidean region into  a very curved one, no matter what happens
far away. The results of sections 1 through 10 require no
dimensional or curvature restrictions, and are not immediately
related to Hamilton program for geometrization of three manifolds.
\par The work on details of this program starts in $\S 11,$ where
we describe the ancient solutions with nonnegative curvature that
may occur as blow-up limits of finite time singularities ( they
must satisfy a certain noncollapsing assumption, which, in the
interpretation of $\S 5,$ corresponds to having bounded entropy).
Then in $\S 12$ we describe the regions of high curvature under
the assumption of almost nonnegative curvature, which is
guaranteed to hold by the Hamilton and Ivey result, mentioned
above. We also prove, under the same assumption, some results on
the control of the curvatures forward and backward in time in
terms of the curvature and volume at a given time in a given ball.
Finally, in $\S 13$ we give a brief sketch of the proof of
geometrization conjecture.
\par The subsections marked by * contain historical remarks and
references. See also [Cao-C] for a relatively recent survey on the
Ricci flow.
\end{intropar}
\section{Ricci flow as a gradient flow}
\begin{parag} Consider the functional $\cF=\int_M{(R+|\nabla
f|^2)e^{-f}dV}$ for a riemannian metric $g_{ij}$ and a function
$f$ on a closed manifold $M$. Its first variation can be expressed
as follows:
$$ \delta \cF (v_{ij},h)=\int_M e^{-f}[-\triangle
v+\nabla_i\nabla_jv_{ij}-R_{ij}v_{ij}$$ $$-v_{ij}\nabla_i
f\nabla_j f+2<\nabla f,\nabla h>+(R+|\nabla f|^2)(v/2-h)]$$
$$=\int_M{e^{-f}[-v_{ij}(R_{ij}+\nabla_i\nabla_j
f)+(v/2-h)(2\triangle f-|\nabla f|^2+R)]},$$where $\delta
g_{ij}=v_{ij}$, $\delta f=h$, $v=g^{ij}v_{ij}$. Notice that
$v/2-h$ vanishes identically iff the measure $dm=e^{-f}dV$ is kept
fixed. Therefore, the symmetric tensor $-(R_{ij}+\nabla_i\nabla_j
f)$ is the $L^2$ gradient of the functional $\cF^m
=\int_M{(R+|\nabla f|^2)dm}$, where now $f$ denotes $\log(dV/dm)$.
Thus given a measure $m$ , we may consider the gradient flow
$(g_{ij})_t=-2(R_{ij}+\nabla_i\nabla_j f)$ for $\cF^m$. For
general $m$ this flow may not exist even for short time; however,
when it exists, it is just the Ricci flow, modified by a
diffeomorphism. The remarkable fact here is that different choices
of $m$ lead to the same flow, up to a diffeomorphism; that is, the
choice of $m$ is analogous to the choice of gauge.
\end{parag}
\begin{prop}
  Suppose that the gradient flow for $\cF^m$ exists for
$t\in[0,T].$ Then at $t=0$ we have $\cF^m\le
\frac{n}{2T}\int_M{dm}.$ \end{prop} {\it Proof.}   We may assume
$\int_M{dm}=1.$ The evolution equations for the gradient flow of
$\cF^m$ are \begin{equation} (g_{ij})_t=-2(R_{ij}+\nabla_i\nabla_j
f) ,\ \  f_t=-R-\triangle f ,\end{equation} and $\cF^m$ satisfies
\begin{equation}
 \cF^m_t=2\int{|R_{ij}+\nabla_i\nabla_j f|^2 dm}
\end{equation}
Modifying by an appropriate diffeomorphism, we get evolution
equations \begin{equation}   (g_{ij})_t=-2R_{ij} , f_t=-\triangle
f + |\nabla f|^2 - R ,\end{equation} and retain (1.2) in the form
\begin{equation} \cF_t=2\int{|R_{ij}+\nabla_i\nabla_j f|^2e^{-f}dV}
\end{equation}  Now we compute
$$\cF_t\ge\frac{2}{n}\int{(R+\triangle
f)^2e^{-f}dV}\ge\frac{2}{n}(\int{(R+\triangle
f)e^{-f}dV})^2=\frac{2}{n}\cF^2,$$ and the proposition follows.

\par {\bf 1.3} {\it Remark.}
 The functional $\cF^m $ has a natural interpretation in terms of
Bochner-Lichnerovicz formulas. The classical formulas of Bochner
(for one-forms) and Lichnerovicz (for spinors) are $\nabla^*\nabla
u_i=(d^*d+dd^*)u_i-R_{ij}u_j$ and $\nabla^*\nabla
\psi=\delta^2\psi-1/4 R\psi.$  Here the operators $\nabla^*$ ,
$d^*$ are defined using the riemannian volume form; this volume
form is also implicitly used in the definition of the Dirac
operator $\delta$ via the requirement $\delta^*=\delta.$  A
routine computation shows that  if we substitute $dm=e^{-f}dV$ for
$dV$ , we get modified Bochner-Lichnerovicz formulas
$\nabla^{*m}\nabla u_i=(d^{*m}d+dd^{*m})u_i-R_{ij}^m u_j $ and
$\nabla^{*m}\nabla\psi=(\delta^m)^2\psi-1/4R^m\psi,$ where
$\delta^m\psi=\delta\psi-1/2(\nabla f)\cdot\psi$ ,
$R_{ij}^m=R_{ij}+\nabla_i\nabla_j f$ , $R^m=2\triangle f-|\nabla
f|^2 +R .$ Note that $g^{ij}R_{ij}^m= R + \triangle f \ne R^m .$
However, we do have the Bianchi identity
$\nabla_i^{*m}R_{ij}^m=\nabla_iR_{ij}^m-R_{ij}\nabla_i
f=1/2\nabla_jR^m .$ Now
$\cF^m=\int_M{{R^m}dm}=\int_M{g^{ij}R_{ij}^m dm}.$

\par {\bf 1.4*}
 The Ricci flow modified by a diffeomorphism was considered by
DeTurck, who observed that by an appropriate choice of
diffeomorphism one can turn the equation from weakly parabolic
into strongly parabolic, thus considerably simplifying the proof
of short time existence and uniqueness; a nice version of DeTurck
trick can be found in [H 4,$\S 6$].
\par The functional $\cF$ and its first variation formula can be
found in the literature on the string theory, where it describes
the low energy effective action; the function $f$ is called
dilaton field; see [D,$\S 6$] for instance.
\par The Ricci tensor $R_{ij}^m$ for a riemannian manifold with a
smooth measure has been used  by Bakry and Emery [B-Em]. See also
a very recent paper [Lott].

\section{No breathers theorem I} \begin{parag}\par  A metric $g_{ij}(t)$
evolving by the Ricci flow is called a breather, if for some
$t_1<t_2 $ and $\alpha>0$ the metrics $\alpha g_{ij}(t_1)$ and
$g_{ij}(t_2)$ differ only by a diffeomorphism; the cases $\alpha=1
, \alpha<1 , \alpha>1 $ correspond to steady, shrinking and
expanding breathers, respectively. Trivial breathers, for which
the metrics $g_{ij}(t_1)$ and $g_{ij}(t_2)$ differ only by
diffeomorphism and scaling for each pair of $t_1$ and $t_2$, are
called Ricci solitons. (Thus, if one considers Ricci flow as a
dynamical system on the space of riemannian metrics modulo
diffeomorphism and scaling, then breathers and solitons correspond
to periodic orbits and fixed points respectively). At each time
the Ricci soliton metric satisfies an equation of the form
$R_{ij}+cg_{ij}+\nabla_i b_j +\nabla_j b_i=0,$ where $c$ is a
number and $b_i$ is a one-form; in particular, when
$b_i=\frac{1}{2}\nabla_i a$ for some function $a$ on $M,$ we get a
gradient Ricci soliton. An important example of a gradient
shrinking soliton is the Gaussian soliton, for which the metric
$g_{ij}$ is just the euclidean metric on $\mathbb{R}^n$, $c=1$ and
$a=-|x|^2/2.$
\par In this and the next section we use the gradient
interpretation of the Ricci flow to rule out nontrivial breathers
(on closed $M$). The argument in the steady case is pretty
straightforward; the expanding case is a little bit more subtle,
because our functional $\cF$ is not scale invariant. The more
difficult shrinking case is discussed in section
3.\end{parag}\begin{parag}
\par  Define $\lambda(g_{ij})=\mbox{inf}\  \cF(g_{ij},f) ,$ where
infimum is taken over all smooth $f,$ satisfying $
\int_M{e^{-f}dV}=1 .$ Clearly, $\lambda(g_{ij})$ is just the
lowest eigenvalue of the operator $-4\triangle+R.$ Then formula
(1.4) implies that $\lambda(g_{ij}(t))$ is nondecreasing in $t,$
and moreover, if $\lambda(t_1)=\lambda(t_2),$ then for $t\in
[t_1,t_2]$ we have $R_{ij}+\nabla_i\nabla_j f=0$ for $f$ which
minimizes $\cF.$ Thus a steady breather is necessarily a steady
soliton. \end{parag}\begin{parag}\par  To deal with the expanding
case consider a scale invariant version
$\bar{\lambda}(g_{ij})=\lambda(g_{ij})V^{2/n}(g_{ij}).$ The
nontrivial expanding breathers will be ruled out once we prove the
following
\par {\bf Claim} {\it $\bar{\lambda}$ is nondecreasing along the Ricci flow
whenever it is nonpositive; moreover, the monotonicity is strict
unless we are on a gradient soliton.} \par (Indeed, on an
expanding breather we would necessarily have $dV/dt>0$ for some $t
{\in } [t_1,t_2].$ On the other hand, for every $t$,
$-\frac{d}{dt}\mbox{log}V=\frac{1}{V}\int{RdV}\ge\lambda(t),$ so
$\bar{\lambda}$ can not be nonnegative everywhere on $[t_1,t_2], $
and the claim applies.)
\par {\it Proof of the claim.}
$${\small\begin{array}{cc}d\bar{\lambda}(t)/dt\ge2V^{2/n}\int{|R_{ij}+\nabla_i\nabla_j
f|^2e^{-f}dV}+\frac{2}{n}V^{(2-n)/n}\lambda\int{-RdV}\ge\\\\2V^{2/n}[\int{|R_{ij}+\nabla_i\nabla_j
f-\frac{1}{n}(R+\triangle
f)g_{ij}|^2e^{-f}dV}+\\\\\frac{1}{n}(\int{(R+\triangle
f)^2e^{-f}dV}-(\int{(R+\triangle
f)e^{-f}dV})^2)]\ge0,\end{array}}$$ where $f$ is the minimizer for
$\cF.$ \end{parag}\begin{parag}\par  The arguments above also show
that there are no nontrivial (that is with non-constant Ricci
curvature) steady or expanding Ricci solitons (on closed $M$).
Indeed, the equality case in the chain of inequalities above
requires that $R+\triangle f$ be constant on $M$; on the other
hand, the Euler-Lagrange equation for the minimizer $f$ is
$2\triangle f-|\nabla f|^2+R=const.$ Thus, $\triangle f-|\nabla
f|^2=const=0$, because $\int{(\triangle f-|\nabla
f|^2)e^{-f}dV}=0.$ Therefore, $f$ is constant by the maximum
principle.\end{parag}\begin{parstr}
\par  A similar, but simpler proof of the results in this
section, follows immediately from  [H 6,$\S 2$], where Hamilton
checks that the minimum of $RV^{\frac{2}{n}}$ is nondecreasing
whenever it is nonpositive, and monotonicity is strict unless the
metric has constant Ricci curvature.
\end{parstr}
\section{ No breathers theorem II}
\begin{parag}
\par  In order to handle the shrinking case when $\lambda>0,$
we need to replace our functional $\cF$ by its generalization,
which contains explicit insertions of the scale parameter, to be
denoted by $\tau.$ Thus consider the functional \begin{equation}
\cW(g_{ij},f,\tau)=\int_M{[\tau(|\nabla
f|^2+R)+f-n](4\pi\tau)^{-\frac{n}{2}}e^{-f}dV},\end{equation}
restricted to $f$ satisfying \begin{equation}
\int_M{(4\pi\tau)^{-\frac{n}{2}}e^{-f}dV}=1,\end{equation}
$\tau>0.$ Clearly $\cW$ is invariant under simultaneous scaling of
$\tau$ and $g_{ij}.$ The evolution equations, generalizing (1.3)
are \begin{equation} (g_{ij})_t=-2R_{ij} , f_t=-\triangle
f+|\nabla f|^2-R+\frac{n}{2\tau} , \tau_t=-1 \end{equation} The
evolution equation for $f$ can also be written as follows:
$\Box^*u=0,$ where $u=(4\pi\tau)^{-\frac{n}{2}}e^{-f},$ and
$\Box^*=-\partial/\partial t -\triangle +R$ is the conjugate heat
operator. Now a routine computation gives \begin{equation}
d\cW/dt=\int_M{2\tau|R_{ij}+\nabla_i\nabla_j
f-\frac{1}{2\tau}g_{ij}|^2(4\pi\tau)^{-\frac{n}{2}}e^{-f}dV}
.\end{equation} Therefore, if we let $\mu(g_{ij},\tau)=\mbox{inf}\
\cW(g_{ij},f,\tau)$ over smooth $f$ satisfying (3.2), and
$\nu(g_{ij})=\mbox{inf}\  \mu(g_{ij},\tau) $ over all positive
$\tau,$ then $\nu(g_{ij}(t))$ is nondecreasing along the Ricci
flow. It is not hard to show that in the definition of $\mu$ there
always exists a smooth minimizer $f$ (on a closed $M$). It is also
clear that $\lim_{\tau\to\infty}\mu(g_{ij},\tau)=+\infty$ whenever
the first eigenvalue of $-4\triangle +R$ is positive. Thus, our
statement that there is no shrinking breathers other than gradient
solitons, is implied by the following \par {\bf Claim} {\it  For
an arbitrary metric $g_{ij}$ on a closed manifold M, the function
$\mu(g_{ij},\tau)$ is negative for small $\tau>0$ and tends to
zero as $\tau$ tends to zero.} \par {\it Proof of the Claim.
(sketch)} Assume that $\bar{\tau}>0$ is so small that Ricci flow
starting from $g_{ij}$ exists on $[0,\bar{\tau}].$ Let
$u=(4\pi\tau)^{-\frac{n}{2}}e^{-f}$ be the solution of the
conjugate heat equation, starting from a $\delta$-function at
$t=\bar{\tau}, \tau(t)=\bar{\tau}-t.$ Then
$\cW(g_{ij}(t),f(t),\tau(t))$ tends to zero as $t$ tends to
$\bar{\tau},$ and therefore
$\mu(g_{ij},\bar{\tau})\le\cW(g_{ij}(0),f(0),\tau(0))<0 $ by
(3.4).\par Now let $\tau\to 0$ and assume that $f^{\tau}$ are the
minimizers, such that
$$\cW(\frac{1}{2}\tau^{-1}g_{ij},f^{\tau},\frac{1}{2})
=\cW(g_{ij},f^{\tau},\tau)=\mu(g_{ij},\tau)\le c<0.$$ The metrics
$\frac{1}{2}\tau^{-1}g_{ij} $ "converge" to the euclidean metric,
and if we could extract a converging subsequence from $f^{\tau},$
we would get a function $f$ on $\mathbb{R}^n$, such that
$\int_{\mathbb{R}^n}{(2\pi)^{-\frac{n}{2}}e^{-f}dx}=1$ and $$
\int_{\mathbb{R}^n}{[\frac{1}{2}|\nabla
f|^2+f-n](2\pi)^{-\frac{n}{2}}e^{-f}dx}<0 $$ The latter inequality
contradicts the Gaussian logarithmic Sobolev inequality, due to
L.Gross. (To pass to its standard form, take $f=|x|^2/2-2\log\phi$
and integrate by parts) This argument is not hard to make
rigorous; the details are left  to the reader.
\end{parag}

\par {\bf 3.2} {\it  Remark.}  Our monotonicity formula (3.4) can in fact be
used to prove a version of the logarithmic Sobolev inequality
(with description of the equality cases) on shrinking Ricci
solitons. Indeed, assume that a metric $g_{ij}$ satisfies
$R_{ij}-g_{ij}-\nabla_i b_j-\nabla_j b_i=0.$ Then under Ricci
flow, $g_{ij}(t)$ is isometric to $(1-2t)g_{ij}(0),$\ \  $
\mu(g_{ij}(t),\frac{1}{2}-t)=\mu(g_{ij}(0),\frac{1}{2}),$ and
therefore the monotonicity formula (3.4) implies that the
minimizer $f$ for $\mu(g_{ij},\frac{1}{2})$ satisfies
$R_{ij}+\nabla_i\nabla_j f-g_{ij}=0.$ Of course, this argument
requires the existence of minimizer, and  justification of the
integration by parts; this is easy if $M$ is closed, but can also
be done with more efforts on some complete $M$, for instance when
$M$ is the Gaussian soliton. \par {\bf 3.3*} The no breathers
theorem in dimension three was proved by Ivey [I]; in fact, he
also ruled out nontrivial Ricci solitons; his proof uses the
almost nonnegative curvature estimate, mentioned in the
introduction.
\par Logarithmic Sobolev inequalities is a vast area of research;
see [G] for a survey and bibliography up to the year 1992; the
influence of the curvature was discussed by Bakry-Emery [B-Em]. In
the context of geometric evolution equations, the logarithmic
Sobolev inequality occurs in Ecker [E 1].

\section {No local collapsing theorem I} In this section we
present an application of the monotonicity formula (3.4) to the
analysis of singularities of the Ricci flow.
\par \begin{parag} Let $g_{ij}(t)$ be a smooth solution to the Ricci flow
$(g_{ij})_t=-2R_{ij}$ on $[0,T).$ We say that $g_{ij}(t)$ is
locally collapsing at $T,$ if there is a sequence of times $t_k\to
T$ and a sequence of metric balls $B_k=B(p_k,r_k)$ at times $t_k,$
such that $r_k^2/t_k$ is bounded, $|Rm|(g_{ij}(t_k))\le r_k^{-2}$
in $B_k$ and $r_k^{-n}Vol(B_k)\to 0.$ \par {\bf Theorem.} {\it  If
$M$ is closed and $T<\infty,$ then $g_{ij}(t)$ is not locally
collapsing at $T.$}
\par {\it Proof.}   Assume that there is a sequence of collapsing balls
$B_k=B(p_k,r_k)$ at times $t_k\to T.$ Then we claim that
$\mu(g_{ij}(t_k),r_k^2)\to -\infty.$ Indeed one can take
$f_k(x)=-\log\phi(\mbox{dist}_{t_k}(x,p_k)r_k^{-1})+c_k,$ where
$\phi$ is a function of one variable, equal 1 on $[0,1/2],$
decreasing on $[1/2,1],$ and very close to 0 on $[1,\infty),$ and
$c_k$ is a constant; clearly $c_k\to -\infty$ as
$r_k^{-n}Vol(B_k)\to 0.$ Therefore, applying the monotonicity
formula (3.4), we get $\mu(g_{ij}(0),t_k+r_k^2)\to -\infty.$
However this is impossible, since $t_k+r_k^2$ is bounded.
\end{parag}
\par \begin{parag} {\bf Definition}  {\it We say that a metric $g_{ij}$ is
$\kappa$-noncollapsed on the scale $\rho,$ if every metric ball
$B$ of radius $r<\rho,$ which satisfies $|Rm|(x)\le r^{-2}$ for
every $x\in B,$ has volume at least $\kappa r^n.$} \par It is
clear that a limit of $\kappa$-noncollapsed metrics on the scale
$\rho$ is also $\kappa$-noncollapsed on the scale $\rho;$ it is
also clear that $\alpha^2g_{ij}$ is $\kappa$-noncollapsed on the
scale $\alpha\rho$ whenever $g_{ij}$ is $\kappa$-noncollapsed on
the scale $\rho.$ The theorem above essentially says that given a
metric $g_{ij}$ on a closed manifold $M$ and $T<\infty,$ one can
find $\kappa=\kappa(g_{ij},T)>0,$ such that the solution
$g_{ij}(t)$ to the Ricci flow starting at $g_{ij}$ is
$\kappa$-noncollapsed on the scale $T^{1/2}$ for all $t\in [0,T),$
provided it exists on this interval. Therefore, using the
convergence theorem of Hamilton, we obtain the following
\par {\bf Corollary.}  {\it Let $g_{ij}(t), t\in [0,T)$ be a solution to the
Ricci flow on a closed manifold $M,$ $T<\infty.$ Assume that for
some sequences $t_k\to T, p_k\in M$ and some constant $C$  we have
$Q_k=|Rm|(p_k,t_k)\to\infty$ and $ |Rm|(x,t)\le CQ_k,$ whenever
$t<t_k.$ Then (a subsequence of) the scalings of $g_{ij}(t_k)$ at
$p_k$ with factors $Q_k$ converges to a complete ancient solution
to the Ricci flow, which is $\kappa$-noncollapsed on all scales
for some $\kappa>0.$}  \end{parag}  \section {A statistical
analogy} In this section we show that the functional $\cW,$
introduced in section 3, is in a sense analogous to minus entropy.
\par {\bf 5.1} Recall that the partition function for the canonical
ensemble at temperature $\beta^{-1}$ is given by
$Z=\int{exp(-\beta E)d\omega(E)},$ where $\omega(E)$ is a "density
of states" measure, which does not depend on $\beta.$ Then one
computes the average energy
$<E>=-\frac{\partial}{\partial\beta}\log Z,$ the entropy
$S=\beta<E>+\log Z,$ and the fluctuation
$\sigma=<(E-<E>)^2>=\frac{\partial^2}{(\partial\beta)^2}\log Z.$
\par Now fix a closed manifold $M$ with a probability measure $m$,
and suppose that our system is described by a metric
$g_{ij}(\tau),$ which depends on the temperature $\tau$ according
to equation $(g_{ij})_\tau=2(R_{ij}+\nabla_i\nabla_j f),$ where
$dm=udV, u=(4\pi\tau)^{-\frac{n}{2}}e^{-f},$ and the partition
function is given by $\log Z=\int{(-f+\frac{n}{2})dm}.$ (We do not
discuss here what assumptions on $g_{ij}$ guarantee that the
corresponding "density of states" measure can be found) Then we
compute $$ <E>=-\tau^2\int_M{(R+|\nabla
f|^2-\frac{n}{2\tau})dm},$$
$$S=-\int_M{(\tau(R+|\nabla f|^2)+f-n)dm},$$
$$ \sigma=2\tau^4\int_M{|R_{ij}+\nabla_i\nabla_j
f-\frac{1}{2\tau}g_{ij}|^2dm}$$ \par Alternatively, we could
prescribe the evolution equations by replacing the $t$-derivatives
by minus $\tau$-derivatives in (3.3 ), and get the same formulas
for $Z, <E>, S, \sigma,$ with $dm$ replaced by $udV.$ \par
Clearly, $\sigma$ is nonnegative; it vanishes only on a gradient
shrinking soliton. $<E>$ is nonnegative as well, whenever the flow
exists for all sufficiently small $\tau>0$ (by proposition 1.2).
Furthermore, if (a) $u$ tends to a $\delta$-function as $\tau\to
0,$ or (b) $u$ is a limit of a sequence of functions $u_i,$ such
that each $u_i$ tends to a $\delta$-function as $\tau\to\tau_i>0,$
and $\tau_i\to 0,$ then $S$ is also nonnegative. In case (a) all
the quantities $<E>, S, \sigma$ tend to zero as $\tau\to 0,$ while
in case (b), which may be interesting if $g_{ij}(\tau)$ goes
singular at $\tau=0,$ the entropy $S$ may tend to a positive
limit. \par If the flow is defined for all sufficiently large
$\tau$ (that is, we have an ancient solution to the Ricci flow, in
Hamilton's terminology), we may be interested in the behavior of
the entropy $S$ as $\tau\to\infty.$  A natural question is whether
we have a gradient shrinking soliton whenever $S$ stays bounded.
\par {\bf 5.2} {\it Remark.}  Heuristically, this statistical analogy is related to
the description of the renormalization group flow, mentioned in
the introduction: in the latter one obtains various quantities by
averaging over higher energy states, whereas in the former those
states are  suppressed by the exponential factor.
\par {\bf 5.3*} An entropy formula for the Ricci flow in dimension two
was found by Chow [C]; there seems to be no relation between his
formula and ours.
\par The interplay of statistical physics and (pseudo)-riemannian
geometry occurs in the subject of Black Hole Thermodynamics,
developed by Hawking et al. Unfortunately, this subject is beyond
my understanding at the moment.
\section { Riemannian formalism in potentially infinite
dimensions} When one is talking of the canonical ensemble, one is
usually considering an embedding of the system of interest into a
much larger standard system of fixed temperature (thermostat). In
this section we attempt to describe such an embedding using the
formalism of Rimannian geometry. \par {\bf 6.1}  Consider the
manifold $\tilde{M}=M\times\mathbb{S}^N\times\mathbb{R}^+$ with
the following metric:$$ \tilde{g}_{ij}=g_{ij},
\tilde{g}_{\alpha\beta}=\tau g_{\alpha\beta},
\tilde{g}_{00}=\frac{N}{2\tau}+R, \tilde{g}_{i\alpha}=\tilde{g}_{i
0}=\tilde{g}_{\alpha 0}=0,$$ where $i,j$ denote coordinate indices
on the $M$ factor, $\alpha,\beta$ denote those on the
$\mathbb{S}^N$ factor, and the coordinate $\tau$ on $\mathbb{R}^+$
has index $0$; $g_{ij}$ evolves with $\tau$ by the backward Ricci
flow $(g_{ij})_\tau=2R_{ij},$ $g_{\alpha\beta}$ is the metric on
$\mathbb{S}^N$ of constant curvature $\frac{1}{2N}.$ It turns out
that the components of the curvature tensor of this metric
coincide (modulo $N^{-1}$) with the components of the matrix
Harnack expression (and its traces), discovered by Hamilton [H 3].
One can also compute that all the components of the Ricci tensor
are equal to zero (mod $N^{-1}$). The heat equation and the
conjugate heat equation on $M$ can be interpreted via Laplace
equation on $\tilde{M}$ for functions and volume forms
respectively: $u$ satisfies the heat equation on $M$ iff
$\tilde{u}$ (the extension of $u$ to $\tilde{M}$ constant along
the $\mathbb{S}^N$ fibres) satisfies
$\tilde{\triangle}\tilde{u}=0\ \mbox{mod}\ N^{-1};$ similarly, $u$
satisfies the conjugate heat equation on $M$ iff
$\tilde{u}^*=\tau^{-\frac{N-1}{2}}\tilde{u}$ satisfies
$\tilde{\triangle}\tilde{u}^*=0\ \ \mbox{mod}\ N^{-1}$ on
$\tilde{M}.$ \par {\bf 6.2}  Starting from $\tilde{g},$ we can
also construct a metric $g^m$ on $\tilde{M},$ isometric to
$\tilde{g}$ (mod $N^{-1}$), which corresponds to the backward
$m$-preserving Ricci flow ( given by equations (1.1) with
$t$-derivatives replaced by minus $\tau$-derivatives,
$dm=(4\pi\tau)^{-\frac{n}{2}}e^{-f}dV$). To achieve this, first
apply to $\tilde{g}$ a (small) diffeomorphism, mapping each point
$(x^{i},y^{\alpha},\tau)$ into $
(x^{i},y^{\alpha},\tau(1-\frac{2f}{N}));$ we would get a metric
$\tilde{g}^m,$ with components (mod $N^{-1}$) $$
\tilde{g}^m_{ij}=\tilde{g}_{ij},
\tilde{g}^m_{\alpha\beta}=(1-\frac{2f}{N})\tilde{g}_{\alpha\beta},
\tilde{g}^m_{00}=\tilde{g}_{00}-2f_{\tau}-\frac{f}{\tau},
\tilde{g}^m_{i 0}=-\nabla_i f, \tilde{g}^m_{i
\alpha}=\tilde{g}^m_{\alpha 0}=0;$$ then apply a horizontal (that
is, along the $M$ factor) diffeomorphism to get $g^m$ satisfying
$(g^m_{ij})_\tau=2(R_{ij}+\nabla_i\nabla_j f);$ the other
components of $g^m$ become (mod $N^{-1}$)
$$g^m_{\alpha\beta}=(1-\frac{2f}{N})\tilde{g}_{\alpha\beta},
g^m_{00}=\tilde{g}^m_{00}-|\nabla
f|^2=\frac{1}{\tau}(\frac{N}{2}-[\tau(2\triangle f-|\nabla f|^2
+R)+f-n]),$$ $$ g^m_{i 0}=g^m_{\alpha 0}=g^m_{i \alpha}=0$$ Note
that the hypersurface $\tau=$const in the metric $g^m$ has the
volume form $\tau^{N/2}e^{-f}$ times the canonical form on $M$ and
$\mathbb{S}^N,$ and the scalar curvature of this hypersurface is
$\frac{1}{\tau}(\frac{N}{2}+\tau(2\triangle f-|\nabla f|^2+R)+f)$
mod $N^{-1}.$  Thus the entropy $S$ multiplied by  the inverse
temperature $\beta$ is essentially minus the total scalar
curvature of this hypersurface. \par {\bf 6.3} Now we return to
the metric $\tilde{g}$ and try to use its Ricci-flatness by
interpreting the Bishop-Gromov relative volume comparison theorem.
Consider a metric ball in $(\tilde{M},\tilde{g})$ centered at some
point $p$ where $\tau=0.$ Then clearly the shortest geodesic
between $p$ and an arbitrary point $q$ is always orthogonal to the
$\mathbb{S}^N$ fibre. The length of such curve $\gamma(\tau)$ can
be computed as
$$\int_0^{\tau(q)}{\sqrt{\frac{N}{2\tau}+R+|\dot{\gamma}_M(\tau)|^2}d\tau}$$
$$
=\sqrt{2N\tau(q)}+\frac{1}{\sqrt{2N}}\int_0^{\tau(q)}{\sqrt{\tau}(R+|\dot{\gamma}_M(\tau)|^2)d\tau}+
O(N^{-\frac{3}{2}})$$ Thus a shortest geodesic should minimize
$\mathcal{L}(\gamma)=\int_0^{\tau(q)}{\sqrt{\tau}(R+|\dot{\gamma}_M(\tau)|^2)d\tau},$
an expression defined entirely in terms of $M$. Let $L(q_M)$
denote the corresponding infimum. It follows that a metric sphere
in $\tilde{M}$ of radius $\sqrt{2N\tau(q)}$ centered at $p$ is
$O(N^{-1})$-close to the hypersurface $\tau=\tau(q),$ and its
volume can be computed as
$V(\mathbb{S}^N)\int_M{(\sqrt{\tau(q)}-\frac{1}{2N}L(x)+O(N^{-2}))^Ndx},$
so the ratio of this volume to $\sqrt{2N\tau(q)}^{N+n}$ is just
constant times $N^{-\frac{n}{2}}$ times
$$\int_M{\tau(q)^{-\frac{n}{2}}\mbox{exp}(-\frac{1}{\sqrt{2\tau(q)}}L(x))dx}+O(N^{-1})$$
The computation suggests that this integral, which we will call
the reduced volume and denote by $\tilde{V}(\tau(q)),$ should be
increasing as $\tau$ decreases. A rigorous proof of this
monotonicity is given in the next section. \par {\bf 6.4*} The
first geometric interpretation of Hamilton's Harnack expressions
was found by Chow and Chu [C-Chu 1,2]; they construct a
potentially degenerate riemannian metric on $M\times \mathbb{R},$
which potentially satisfies the Ricci soliton equation; our
construction is, in a certain sense, dual to theirs. \par Our
formula for the reduced volume resembles the expression in Huisken
monotonicity formula for the mean curvature flow [Hu]; however, in
our case the monotonicity is in the opposite direction.
\section { A comparison geometry approach to the Ricci flow} {\bf 7.1} In
this section we consider an evolving metric
$(g_{ij})_{\tau}=2R_{ij}$ on a manifold $M;$ we assume that either
$M$ is closed, or $g_{ij}(\tau)$ are complete and have uniformly
bounded curvatures. To each curve $\gamma(\tau),
0<\tau_1\le\tau\le\tau_2,$ we associate its $\mathcal{L}$-length
$$\mathcal{L}(\gamma)=\int_{\tau_1}^{\tau_2}{\sqrt{\tau}(R(\gamma(\tau))+|\dot{\gamma}(\tau)|^2)d\tau}$$
(of course, $ R(\gamma(\tau))$ and $|\dot{\gamma}(\tau)|^2$ are
computed using $g_{ij}(\tau)$) \par Let
$X(\tau)=\dot{\gamma}(\tau),$ and let $Y(\tau)$ be any vector
field along $\gamma(\tau).$ Then the first variation formula can
be derived as follows:  $$\delta_Y(\mathcal{L}) =$$
$$\int_{\tau_1}^{\tau_2}{\sqrt{\tau}(<Y,\nabla R>+2<\nabla_Y
X,X>)d\tau} =\int_{\tau_1}^{\tau_2}{\sqrt{\tau}(<Y,\nabla
R>+2<\nabla_X Y,X>)d\tau}$$
$$=\int_{\tau_1}^{\tau_2}{\sqrt{\tau}(<Y,\nabla
R>+2\frac{d}{d\tau}<Y,X>-2<Y,\nabla_X X>-4\mbox{Ric}(Y,X))d\tau}$$
\begin{equation}
=\left.2\sqrt{\tau}<X,Y>\right|_{\tau_1}^{\tau_2}+\int_{\tau_1}^{\tau_2}{\sqrt{\tau}<Y,\nabla
R-2\nabla_X X-4\mbox{Ric}(X,\cdot)-\frac{1}{\tau}X>d\tau}
\end{equation} Thus $\mathcal{L}$-geodesics must satisfy
\begin{equation}
\nabla_X X-\frac{1}{2}\nabla R+
\frac{1}{2\tau}X+2\mbox{Ric}(X,\cdot)=0 \end{equation} Given two
points $p,q$ and $\tau_2>\tau_1>0,$ we can always find an
$\mathcal{L}$-shortest curve $\gamma(\tau),
\tau\in[\tau_1,\tau_2]$ between them, and every such
$\mathcal{L}$-shortest curve is $\mathcal{L}$-geodesic. It is easy
to extend this to the case $\tau_1=0;$ in this case
$\sqrt{\tau}X(\tau)$ has a limit as $\tau\to 0.$ From now on we
fix $p$ and $\tau_1=0$ and denote by $L(q,\bar{\tau})$ the
$\mathcal{L}$-length of the $\mathcal{L}$-shortest curve
$\gamma(\tau), 0\le\tau\le\bar{\tau},$ connecting $p$ and $q.$ In
the computations below we pretend that shortest
$\mathcal{L}$-geodesics between $p$ and $q$ are unique for all
pairs $(q,\bar{\tau}); $ if this is not the case, the inequalities
that we obtain are still valid when understood in the barrier
sense, or in the sense of distributions. \par The first variation
formula (7.1) implies that $\nabla
L(q,\bar{\tau})=2\sqrt{\bar{\tau}}X(\bar{\tau}),$ so that $|\nabla
L|^2=4\bar{\tau}|X|^2=-4\bar{\tau}R+4\bar{\tau}(R+|X|^2).$ We can
also compute
$$L_{\bar{\tau}}(q,\bar{\tau})=\sqrt{\bar{\tau}}(R+|X|^2)-<X,\nabla
L>=2\sqrt{\bar{\tau}}R-\sqrt{\bar{\tau}}(R+|X|^2)$$ To evaluate
$R+|X|^2$ we compute (using (7.2))
$$  \frac{d}{d\tau}(R(\gamma(\tau))+|X(\tau)|^2)=R_{\tau}+<\nabla
R,X>+2<\nabla_X X,X>+2\mbox{Ric}(X,X)$$
$$=R_{\tau}+\frac{1}{\tau}R+2<\nabla
R,X>-2\mbox{Ric}(X,X)-\frac{1}{\tau}(R+|X|^2)$$
\begin{equation}=-H(X)-\frac{1}{\tau}(R+|X|^2),\end{equation} where $H(X)$ is the
Hamilton's expression for the trace Harnack inequality (with
$t=-\tau$). Hence, \begin{equation}
\bar{\tau}^{\frac{3}{2}}(R+|X|^2)(\bar{\tau})=-K+\frac{1}{2}L(q,\bar{\tau}),\end{equation}
where $K=K(\gamma,\bar{\tau}) $ denotes the integral
$\int_0^{\bar{\tau}}{\tau^{\frac{3}{2}}H(X)d\tau},$ which we'll
encounter a few times below. Thus we get \begin{equation}
L_{\bar{\tau}}=2\sqrt{\bar{\tau}}R-\frac{1}{2\bar{\tau}}L+\frac{1}{\bar{\tau}}K
\end{equation}
\begin{equation} |\nabla
L|^2=-4\bar{\tau}R+\frac{2}{\sqrt{\bar{\tau}}}L-\frac{4}{\sqrt{\bar{\tau}}}K
\end{equation}
\par Finally we need to estimate the second variation of $L.$ We
compute
$$\de^2_Y(\mathcal{L})=\int_0^{\bar{\tau}}{\sqrt{\tau}(Y\cdot Y\cdot
R+2<\nabla_Y \nabla_Y X,X>+2|\nabla_Y X|^2)d\tau}$$
$$=\int_0^{\bar{\tau}}{\sqrt{\tau}(Y\cdot Y\cdot R+2<\nabla_X
\nabla_Y Y,X>+2<R(Y,X),Y,X>+2|\nabla_X Y|^2)d\tau}$$ Now $$
\frac{d}{d\tau}<\nabla_Y Y,X>=<\nabla_X \nabla_Y Y,X>+<\nabla_Y
Y,\nabla_X X>+2Y\cdot\mbox{Ric}(Y,X)-X\cdot\mbox{Ric}(Y,Y),$$ so,
if $Y(0)=0$ then
$$ \de^2_Y(\mathcal{L})=2<\nabla_Y
Y,X>\sqrt{\bar{\tau}}+$$
\begin{multline}\int_0^{\bar{\tau}}\sqrt{\tau}(\nabla_Y \nabla_Y
R+2<R(Y,X),Y,X>+2|\nabla_X
Y|^2\\+2\nabla_X\mbox{Ric}(Y,Y)-4\nabla_Y\mbox{Ric}(Y,X))d\tau,\end{multline}
where we discarded the scalar product of $-2\nabla_Y Y$ with the
left hand side of (7.2). Now fix the value of $Y$ at
$\tau=\bar{\tau}$, assuming $|Y(\bar{\tau})|=1,$ and construct $Y$
on $[0,\bar{\tau}] $ by solving the ODE \begin{equation} \nabla_X
Y=-\mbox{Ric}(Y,\cdot)+\frac{1}{2\tau}Y \end{equation} We compute
$$ \frac{d}{d\tau}<Y,Y>=2\mbox{Ric}(Y,Y)+2<\nabla_X
Y,Y>=\frac{1}{\tau}<Y,Y>,$$ so
$|Y(\tau)|^2=\frac{\tau}{\bar{\tau}},$ and in particular,
$Y(0)=0.$ Making a substitution into (7.7), we get
$$\mbox{Hess}_L(Y,Y)\le$$
$$\int_0^{\bar{\tau}}\sqrt{\tau}(\nabla_Y \nabla_Y
R+2<R(Y,X),Y,X>+2\nabla_X\mbox{Ric}(Y,Y)-4\nabla_Y\mbox{Ric}(Y,X)$$
$$
 +2|\mbox{Ric}(Y,\cdot)|^2-
\frac{2}{\tau}\mbox{Ric}(Y,Y)+\frac{1}{2\tau\bar{\tau}})d\tau$$ To
put this in a more convenient form, observe that
$$\frac{d}{d\tau}\mbox{Ric}(Y(\tau),Y(\tau))=\mbox{Ric}_{\tau}(Y,Y)+\nabla_X\mbox{Ric}(Y,Y)+
2\mbox{Ric}(\nabla_X Y,Y)$$
$$=\mbox{Ric}_{\tau}(Y,Y)+\nabla_X\mbox{Ric}(Y,Y)+\frac{1}{\tau}\mbox{Ric}(Y,Y)-
2|\mbox{Ric}(Y,\cdot)|^2,$$ so \begin{equation}
\mbox{Hess}_L(Y,Y)\le\frac{1}{\sqrt{\bar{\tau}}}-2\sqrt{\bar{\tau}}\mbox{Ric}(Y,Y)-\int_0^{\bar{\tau}}
{\sqrt{\tau}H(X,Y)d\tau},\end{equation} where $$H(X,Y)=-\nabla_Y
\nabla_Y
R-2<R(Y,X)Y,X>-4(\nabla_X\mbox{Ric}(Y,Y)-\nabla_Y\mbox{Ric}(Y,X))$$
$$-2\mbox{Ric}_{\tau}(Y,Y)+
2|\mbox{Ric}(Y,\cdot)|^2-\frac{1}{\tau}\mbox{Ric}(Y,Y)$$ is the
Hamilton's expression for the matrix Harnack inequality (with
$t=-\tau$). Thus \begin{equation} \triangle
L\le-2\sqrt{\tau}R+\frac{n}{\sqrt{\tau}}-\frac{1}{\tau}K\end{equation}
\par A field $Y(\tau)$ along $\mathcal{L}$-geodesic $\gamma(\tau)$
is called $\mathcal{L}$-Jacobi, if it is the derivative of a
variation of $\gamma$ among $\mathcal{L}$-geodesics. For an
$\mathcal{L}$-Jacobi field $Y$ with $|Y(\bar{\tau})|=1$ we have
$$ \frac{d}{d\tau}|Y|^2=2\mbox{Ric}(Y,Y)+2<\nabla_X
Y,Y>=2\mbox{Ric}(Y,Y)+2<\nabla_Y X,Y>$$
\begin{equation}=2\mbox{Ric}(Y,Y)+\frac{1}{\sqrt{\bar{\tau}}}\mbox{Hess}_L(Y,Y)\le\frac{1}{\bar{\tau}}-
\frac{1}{\sqrt{\bar{\tau}}}\int_0^{\bar{\tau}}{\tau^{\frac{1}{2}}H(X,\tilde{Y})d\tau},\end{equation}
where $\tilde{Y}$ is obtained by solving ODE (7.8) with initial
data $\tilde{Y}(\bar{\tau})=Y(\bar{\tau}).$ Moreover, the equality
in (7.11) holds only if $\tilde{Y}$ is $\mathcal{L}$-Jacobi and
hence
$\frac{d}{d\tau}|Y|^2=2\mbox{Ric}(Y,Y)+\frac{1}{\sqrt{\bar{\tau}}}\mbox{Hess}_L(Y,Y)
=\frac{1}{\bar{\tau}}.$ \par Now we can deduce an estimate for the
jacobian $J$ of the $\mathcal{L}$-exponential map, given by
$\mathcal{L}\mbox{exp}_X(\bar{\tau})=\gamma(\bar{\tau}),$ where
$\gamma(\tau)$ is the $\mathcal{L}$-geodesic, starting at $p$ and
having $X$ as the limit of $\sqrt{\tau}\dot{\gamma}(\tau)$ as
$\tau\to 0.$ We obtain
\begin{equation}
\frac{d}{d\tau}\mbox{log}J(\tau)\le\frac{n}{2\bar{\tau}}-\frac{1}{2}\bar{\tau}^{-\frac{3}{2}}K,\end{equation}
with equality only if
$2\mbox{Ric}+\frac{1}{\sqrt{\bar{\tau}}}\mbox{Hess}_L=\frac{1}{\bar{\tau}}g.$
Let $l(q,\tau)=\frac{1}{2\sqrt{\tau}}L(q,\tau)$ be the reduced
distance. Then along an $\mathcal{L}$-geodesic $\gamma(\tau)$ we
have (by (7.4))
$$\frac{d}{d\tau}l(\tau)=-\frac{1}{2\bar{\tau}}l+\frac{1}{2}(R+|X|^2)
=-\frac{1}{2}\bar{\tau}^{-\frac{3}{2}}K,$$ so (7.12) implies that
$\tau^{-\frac{n}{2}}\mbox{exp}(-l(\tau))J(\tau)$ is nonincreasing
in $\tau$ along $\gamma$, and monotonicity is strict unless we are
on a gradient shrinking soliton. Integrating over $M$, we get
monotonicity of the reduced volume function
$\tilde{V}(\tau)=\int_M{\tau^{-\frac{n}{2}}\mbox{exp}(-l(q,\tau))dq}.$
 ( Alternatively, one could obtain the same monotonicity by
 integrating the differential inequality \begin{equation}
 l_{\bar{\tau}}-\triangle l+|\nabla l|^2-R+\frac{n}{2\bar{\tau}}\ge0,\end{equation} which
 follows immediately from (7.5), (7.6) and (7.10). Note also a useful inequality
  \begin{equation} 2\triangle l-|\nabla l|^2 +R+\frac{l-n}{\bar{\tau}}\le
  0,\end {equation}
  which follows from (7.6), (7.10).) \par On the
 other hand, if we denote $\bar{L}(q,\tau)=2\sqrt{\tau}L(q,\tau),$
 then from (7.5), (7.10) we obtain \begin{equation}
 \bar{L}_{\bar{\tau}}+\triangle \bar{L}\le2n\end{equation} Therefore, the
 minimum of $\bar{L}(\cdot,\bar{\tau})-2n\bar{\tau}$ is
 nonincreasing, so in particular, the minimum of
 $l(\cdot,\bar{\tau})$ does not exceed $\frac{n}{2}$ for each
 $\bar{\tau}>0.$ (The lower bound for $l$ is much easier to obtain
 since the evolution equation $R_{\tau}=-\triangle
 R-2|\mbox{Ric}|^2$ implies
 $R(\cdot,\tau)\ge-\frac{n}{2(\tau_0-\tau)},$ whenever the flow
 exists for $\tau\in[0,\tau_0].$) \par {\bf 7.2} If the metrics
 $g_{ij}(\tau) $ have nonnegative curvature operator, then
 Hamilton's differential Harnack inequalities hold, and one can
 say more about the behavior of $l.$ Indeed, in this case, if the
 solution is defined for $\tau\in[0,\tau_0],$ then $
 H(X,Y)\ge-\mbox{Ric}(Y,Y)(\frac{1}{\tau}+\frac{1}{\tau_0-\tau})\ge
 -R(\frac{1}{\tau}+\frac{1}{\tau_0-\tau})|Y|^2$ and
 $H(X)\ge-R(\frac{1}{\tau}+\frac{1}{\tau_0-\tau}).$ Therefore,
 whenever $\tau$ is bounded away from $\tau_0$ (say,
 $\tau\le(1-c)\tau_0, c>0$), we get (using (7.6), (7.11)) \begin{equation} |\nabla l|^2+R\le\frac{Cl}{\tau},\end{equation} and for
 $\mathcal{L}$-Jacobi fields $Y$ \begin{equation}
 \frac{d}{d\tau}\mbox{log}|Y|^2\le\frac{1}{\tau}(Cl+1)\end{equation} \par {\bf
 7.3}
 As the first application of the comparison inequalities above,
 let us give an alternative proof of a weakened version of the no
 local collapsing theorem 4.1. Namely, rather than assuming
 $|Rm|(x,t_k)\le r_k^{-2}$ for $x\in B_k,$ we require
 $|Rm|(x,t)\le r_k^{-2}$ whenever $x\in B_k, t_k-r_k^2\le t\le
 t_k.$ Then the proof can go as follows: let $\tau_k(t)=t_k-t,
 p=p_k,\epsilon_k=r_k^{-1}Vol(B_k)^{\frac{1}{n}}.$ We claim that
 $\tilde{V}_k(\epsilon_k r_k^2) < 3\epsilon_k^{\frac{n}{2}}$ when
 $k$ is large. Indeed, using the $\mathcal{L}$-exponential map we
 can integrate over $T_pM$ rather than $M;$ the vectors in $T_pM$
 of length at most $\frac{1}{2}\epsilon_k^{-\frac{1}{2}}$ give rise to
 $\mathcal{L}$-geodesics, which can not escape from $B_k$ in time
 $\epsilon_k r_k^2,$ so their contribution to the reduced volume
 does not exceed $2\epsilon_k^{\frac{n}{2}};$ on the other hand,
 the contribution of the longer vectors does not exceed
 $\mbox{exp}(-\frac{1}{2}\epsilon_k^{-\frac{1}{2}})$ by the jacobian
 comparison theorem. However, $\tilde{V}_k(t_k)$ (that is, at
 $t=0$) stays bounded away from zero. Indeed, since
 $\mbox{min}\  l_k(\cdot,t_k-\frac{1}{2}T)\le\frac{n}{2},$ we can
 pick a point $q_k,$ where it is attained, and obtain a universal
 upper bound on $l_k(\cdot, t_k)$ by considering only curves
 $\gamma$ with $\gamma(t_k-\frac{1}{2}T)=q_k,$ and using the fact
 that all geometric quantities in $g_{ij}(t)$ are uniformly bounded when
 $t\in[0,\frac{1}{2}T].$ Since the monotonicity of the reduced
 volume requires $\tilde{V}_k(t_k)\le\tilde{V}_k(\epsilon_k
 r_k^2),$ this is a contradiction. \par A similar argument shows
 that the statement of the corollary in 4.2 can be strengthened by
 adding another property of the ancient solution, obtained as a
 blow-up limit. Namely, we may claim that if, say, this solution
 is defined for $t\in(-\infty,0),$ then for any point $p$ and any
 $t_0>0,$ the reduced volume function $\tilde{V}(\tau),$
 constructed using $p$ and $\tau(t)=t_0-t,$ is bounded below by
 $\kappa.$ \par {\bf 7.4*} The computations in
 this section are just natural modifications
 of those in the classical variational theory of geodesics
 that can be found in any textbook on Riemannian
 geometry; an even closer reference is [L-Y], where
 they use "length", associated to a linear parabolic equation,
 which is pretty much the same as in our case.
  \section {No local collapsing theorem II} {\bf 8.1} Let us first
 formalize the notion of local collapsing, that was used in 7.3.
 \par {\bf Definition.} {\it A solution to the Ricci flow
 $(g_{ij})_t=-2R_{ij}$ is said to be $\kappa$-collapsed at
 $(x_0,t_0)$ on the scale $r>0$ if $|Rm|(x,t)\le r^{-2}$ for all
 $(x,t)$ satisfying
 $\mbox{dist}_{t_0}(x,x_0)<r$ and $t_0-r^2\le t\le t_0,$ and the
 volume of the metric ball $B(x_0,r^2)$ at time $t_0$ is less than
 $\kappa r^n.$} \par {\bf 8.2} {\bf Theorem.} {\it  For any $A>0$ there exists
 $\kappa=\kappa(A)>0$ with the following property. If $g_{ij}(t)$
 is a smooth solution to the Ricci flow $(g_{ij})_t=-2R_{ij}, 0\le
 t\le r_0^2,$ which has $|Rm|(x,t)\le r_0^{-2}$  for all $(x,t),$
 satisfying $\mbox{dist}_0(x,x_0)<r_0,$ and the volume of the
 metric ball $B(x_0,r_0)$ at time zero is at least $A^{-1}r_0^n,$
 then $g_{ij}(t) $ can not be $\kappa$-collapsed on the scales
 less than $r_0$ at a point $(x,r_0^2)$ with
 $\mbox{dist}_{r_0^2}(x,x_0)\le Ar_0.$}  \par {\it Proof.}  By scaling we
 may assume $r_0=1;$ we may also assume $\mbox{dist}_1(x,x_0)=A.$
 Let us apply the constructions of 7.1 choosing $p=x,
 \tau(t)=1-t.$ Arguing as in 7.3, we see that if our solution is
 collapsed at $x$ on the scale $r\le 1,$ then the reduced volume
 $\tilde{V}(r^2)$ must be very small; on the other hand,
 $\tilde{V}(1)$ can not be small unless $\mbox{min}\
  l(x,\frac{1}{2})$ over $x$ satisfying
 $\mbox{dist}_{\frac{1}{2}}(x,x_0)\le\frac{1}{10}$ is large. Thus
 all we need is to estimate $l,$ or equivalently $\bar{L},$ in
 that ball. Recall that $\bar{L}$ satisfies the differential
 inequality (7.15). In order to use it efficiently in a maximum
 principle argument, we need first to check the following simple
 assertion.
 \par {\bf 8.3 Lemma.} {\it Suppose we have a solution to the Ricci flow
 $(g_{ij})_t=-2R_{ij}.$ \par (a) Suppose  $\mbox{Ric}(x,t_0)\le
 (n-1)K$ when $ \mbox{dist}_{t_0}(x,x_0)<r_0.$ Then the distance
 function $d(x,t)=\mbox{dist}_t(x,x_0)$ satisfies at $t=t_0$ outside $B(x_0,r_0)$ the
 differential inequality $$ d_t-\triangle d\ge
 -(n-1)(\frac{2}{3}Kr_0+r_0^{-1}) $$
 (the inequality must be understood in the barrier sense, when
 necessary)\par (b)} (cf. [H 4,$\S 17$]) {\it  Suppose $\mbox{Ric}(x,t_0)\le (n-1)K$ when $\mbox{dist}_{t_0}(x,x_0)<r_0,$
  or $\mbox{dist}_{t_0}(x,x_1)<r_0.$ Then $$\frac{d}{dt}\mbox{dist}_t(x_0,x_1)\ge
  -2(n-1)(\frac{2}{3}Kr_0+r_0^{-1}) \ \mbox{at}\ \ \ t=t_0$$} {\it Proof of Lemma.} (a)  Clearly,
 $d_t(x)=\int_{\gamma}{-\mbox{Ric}(X,X)},$ where $\gamma$ is the
 shortest geodesic between $x$ and $x_0$ and $X$ is its unit
 tangent vector, On the other hand, $\triangle d\le
 \sum_{k=1}^{n-1}{s_{Y_k}''(\gamma)},$ where $Y_k$ are vector
 fields along $\gamma,$ vanishing at $x_0$ and forming an
 orthonormal basis at $x$ when complemented by $X,$ and
 $s_{Y_k}''(\gamma)$ denotes the second variation  along $Y_k$
 of the length of $\gamma.$ Take $Y_k$ to be parallel between $x$
 and $x_1,$ and linear between $x_1$ and $x_0,$ where
 $d(x_1,t_0)=r_0.$ Then $$ \triangle
 d\le\sum_{k=1}^{n-1}s_{Y_k}''(\gamma)=\int_{r_0}^{d(x,t_0)}{-\mbox{Ric}(X,X)ds}+\int_0^{r_0}
 {(\frac{s^2}{r_0^2}(-\mbox{Ric}(X,X))+\frac{n-1}{r_0^2})ds} $$
 $$=\int_{\gamma}{-\mbox{Ric}(X,X)}
 +\int_0^{r_0}{(\mbox{Ric}(X,X)(1-\frac{s^2}{r_0^2})+\frac{n-1}{r_0^2})ds}\le
 d_t+(n-1)(\frac{2}{3}Kr_0+r_0^{-1})$$ The proof of (b) is similar. \par Continuing the proof
 of theorem, apply the maximum principle to the function  $
 h(y,t)=\phi(d(y,t)-A(2t-1))(\bar{L}(y,1-t)+2n+1),$ where
 $d(y,t)=\mbox{dist}_t(x,x_0),$ and $\phi$ is a function of one
 variable, equal $1$ on $(-\infty,\frac{1}{20}),$ and rapidly
 increasing to infinity on $(\frac{1}{20},\frac{1}{10}),$ in such
 a way that \begin{equation} 2(\phi ')^2/\phi-\phi ''\ge (2A+100n)\phi
 '-C(A)\phi,\end{equation} for some constant $C(A)<\infty.$ Note that
 $\bar{L}+2n+1\ge 1$ for $t\ge\frac{1}{2}$  by the remark in the
 very end of 7.1. Clearly, $\mbox{min}\  h(y,1)\le h(x,1)=2n+1.$ On
 the other hand, $\mbox{min}\  h(y,\frac{1}{2}) $ is achieved for
 some $y$ satisfying $d(y,\frac{1}{2})\le \frac{1}{10}.$ Now we
 compute \begin{equation} \Box h=(\bar{L}+2n+1)(-\phi
 ''+(d_t-\triangle d-2A)\phi
 ')-2<\nabla\phi\nabla\bar{L}>+(\bar{L}_t-\triangle\bar{L})\phi
 \end{equation}
 \begin{equation} \nabla h=(\bar{L}+2n+1)\nabla\phi+\phi\nabla\bar{L}\end{equation} At
 a minimum point of $h$ we have $\nabla h=0,$ so (8.2) becomes
 \begin{equation} \Box h=(\bar{L}+2n+1)(-\phi ''+(d_t-\triangle d-2A)\phi
 '+2(\phi ')^2/\phi)+(\bar{L}_t-\triangle\bar{L})\phi \end{equation} Now since
 $d(y,t)\ge\frac{1}{20}$ whenever $\phi '\neq 0,$ and
 since $\mbox{Ric}\le n-1$ in $B(x_0,\frac{1}{20}),$ we can apply
 our lemma (a) to get $d_t-\triangle d\ge-100(n-1)$ on the set where
 $\phi '\neq 0.$ Thus, using (8.1) and (7.15), we get $$ \Box
 h\ge-(\bar{L}+2n+1)C(A)\phi-2n\phi\ge-(2n+C(A))h$$ This implies
 that $\mbox{min}\ h$ can not decrease too fast, and we get the
 required estimate. \section {Differential Harnack inequality for
 solutions of the conjugate heat equation} \begin{prop}  Let
 $g_{ij}(t)$ be a solution to the Ricci flow $(g_{ij})_t=-2R_{ij},
 0\le t\le T,$ and let $u=(4\pi(T-t))^{-\frac{n}{2}}e^{-f}$
 satisfy the conjugate heat equation $\Box^*u=-u_t-\triangle
 u+Ru=0.$ Then $v=[(T-t)(2\triangle f-|\nabla f|^2+R)+f-n]u$
 satisfies \begin{equation} \Box^*v=-2(T-t)|R_{ij}+\nabla_i\nabla_j
 f-\frac{1}{2(T-t)}g_{ij}|^2 \end{equation} \end{prop} {\it Proof.}  Routine computation.

 Clearly, this proposition immediately implies the monotonicity
 formula (3.4); its advantage over (3.4) shows up when one has to
 work locally.
 \begin{cor} Under the same assumptions, on a closed manifold $M$,or whenever the
 application of the maximum principle can be justified, $\mbox{min}\ v/u$ is
 nondecreasing in $t.$ \end{cor}\begin{cor} Under the same
 assumptions, if $u$ tends to a $\delta$-function as $t\to T,$ then
 $v\le 0 $ for all $t<T.$\end{cor} {\it Proof.}  If $h$ satisfies the
 ordinary heat equation $h_t=\triangle h$ with respect to the
 evolving metric $g_{ij}(t),$ then we have
 $\frac{d}{dt}\int{hu}=0$ and $\frac{d}{dt}\int{hv}\ge 0.$ Thus we
 only need to check that for everywhere positive $h$ the limit of
 $\int{hv}$ as $t\to T$ is nonpositive. But it is easy to see,
 that this limit is in fact zero. \begin{cor}  Under
 assumptions of the previous corollary, for any smooth curve
 $\gamma(t)$ in $M$ holds \begin{equation}
 -\frac{d}{dt}f(\gamma(t),t)\le\frac{1}{2}(R(\gamma(t),t)+|\dot{\gamma}(t)|^2)
 -\frac{1}{2(T-t)}f(\gamma(t),t) \end{equation} \end{cor}{\it Proof.}  From the evolution
 equation $f_t=-\triangle f+|\nabla f|^2-R+\frac{n}{2(T-t)}$ and
 $v\le 0$ we get $f_t+\frac{1}{2}R-\frac{1}{2}|\nabla
 f|^2-\frac{f}{2(T-t)}\ge 0.$ On the other hand,$
 -\frac{d}{dt}f(\gamma(t),t)=-f_t-<\nabla f,\dot{\gamma}(t)>\le
 -f_t+\frac{1}{2}|\nabla f|^2+\frac{1}{2}|\dot{\gamma}|^2.$
 Summing these two inequalities, we get (9.2).
 \begin{cor}
 If under assumptions of the previous corollary, $p$ is the point
 where the limit $\delta$-function is concentrated, then
 $f(q,t)\le l(q,T-t),$ where $l$ is the reduced distance, defined
 in 7.1, using $p$ and $\tau(t)=T-t.$ \end{cor}{\it Proof.} Use (7.13) in
 the form $\Box^*\mbox{exp}(-l)\le 0.$ \par {\bf 9.6} {\it Remark.} Ricci flow can be
 characterized among all other evolution equations by the
 infinitesimal behavior of the fundamental solutions of the
 conjugate heat equation. Namely, suppose we have a riemannian
 metric $g_{ij}(t)$ evolving with time according to an equation
 $(g_{ij})_t=A_{ij}(t).$ Then we have the heat operator
 $\Box=\frac{\partial}{\partial t}-\triangle$ and its conjugate
 $\Box^*=-\frac{\partial}{\partial t}-\triangle-\frac{1}{2}A,$ so
 that $\frac{d}{dt}\int{uv}=\int{((\Box u)v-u(\Box^* v))}.$
 (Here $A=g^{ij}A_{ij}$) Consider the fundamental solution
 $u=(-4\pi t)^{-\frac{n}{2}}e^{-f}$ for $\Box^*,$ starting as
 $\delta$-function at some point $(p,0).$ Then for general $A_{ij}$
 the function $(\Box\bar{ f}+\frac{\bar{f}}{t})(q,t),$ where $\bar{f}=f-\int{fu},$ is of the order
 $O(1)$ for $(q,t)$ near $(p,0).$ The Ricci flow $A_{ij}=-2R_{ij}$ is
 characterized by the condition $(\Box\bar{
 f}+\frac{\bar{f}}{t})(q,t)=o(1);$ in fact, it is
 $O(|pq|^2+|t|)$ in this case.\par {\bf 9.7*}
 Inequalities of the type of (9.2) are known as
 differential Harnack inequalities; such inequality was proved by Li and Yau [L-Y]
  for the solutions of linear parabolic equations on riemannian manifolds. Hamilton [H 7,8] used
   differential Harnack inequalities for the solutions of backward heat equation on
   a manifold to prove monotonicity formulas for certain parabolic flows. A
  local
  monotonicity formula for mean curvature flow making use of
  solutions of backward heat equation was obtained by Ecker [E 2].

   \section { Pseudolocality theorem}
 \begin{thm}  For every $\alpha>0$ there exist
 $\delta>0,\epsilon>0$ with the following property. Suppose we
 have a smooth solution to the Ricci flow $(g_{ij})_t=-2R_{ij},
 0\le t\le (\epsilon r_0)^2,$ and assume that at $t=0$ we have
 $R(x)\ge -r_0^{-2}$ and
 $Vol(\partial\Omega)^n\ge(1-\delta)c_nVol(\Omega)^{n-1}$ for any
 $x,\Omega\subset B(x_0,r_0),$ where $c_n$ is the euclidean
 isoperimetric constant. Then we have an estimate $|Rm|(x,t)\le
 \alpha t^{-1}+(\e r_0)^{-2}$ whenever $0<t\le (\epsilon
 r_0)^2, d(x,t)=\mbox{dist}_t(x,x_0)<\epsilon r_0.$ \end{thm} Thus,
 under the Ricci flow, the almost singular regions (where
 curvature is large) can not instantly significantly influence the
 almost euclidean regions. Or , using the interpretation via
 renormalization group flow, if a region looks trivial (almost
 euclidean) on higher energy scale, then it can not suddenly
 become highly nontrivial on a slightly lower energy scale.
 \par {\it Proof. }
 It is an argument by contradiction. The idea is to pick a
 point $(\bar{x},\bar{t})$ not far from $(x_0,0)$ and consider the
 solution $u$ to the conjugate heat equation, starting as
 $\delta$-function at $(\bar{x},\bar{t}),$ and the corresponding
 nonpositive function $v$ as in 9.3. If the curvatures at
 $(\bar{x},\bar{t})$ are not small compared to $\bar{t}^{-1}$ and
 are larger than at nearby points, then one can show that
 $\int{v}$  at time $t$ is bounded away from zero for (small) time intervals
 $\bar{t}-t$ of the order of $|Rm|^{-1}(\bar{x},\bar{t}).$ By
 monotonicity we conclude that $\int{v}$ is bounded away from zero
 at $t=0.$ In fact, using (9.1) and an appropriate cut-off
 function, we can show that at $t=0$ already the integral of $v$
 over $B(x_0,r)$ is bounded away from zero, whereas the integral
 of $u$ over this ball is close to $1,$ where $r$ can be made as
 small as we like compared to $r_0.$ Now using the control over
 the scalar curvature and isoperimetric constant in $B(x_0r_0),$
 we can obtain a contradiction to the logarithmic Sobolev
 inequality. \par Now let us go into details. By scaling assume
 that $r_0=1.$ We may also assume that $\alpha$ is small, say
 $\alpha<\frac{1}{100n}.$ From now on we fix $\alpha$ and denote
 by $M_{\alpha}$ the set of pairs $(x,t),$ such that
 $|Rm|(x,t)\ge\alpha t^{-1}.$ \par {\bf Claim 1.}{\it  For any $A>0,$ if
 $g_{ij}(t)$ solves the Ricci flow equation on $0\le t\le
 \epsilon^2, A\epsilon<\frac{1}{100n},$ and $|Rm|(x,t)>\alpha
 t^{-1}+\epsilon^{-2}$ for some $(x,t),$ satisfying $0\le t\le
 \epsilon^2, d(x,t)<\epsilon,$ then one can find $(\bx,\bt)\in
 M_{\al},$ with $0<\bt\le \e^2, d(\bx,\bt)<(2A+1)\e,$ such that \begin{equation}
  |Rm|(x,t)\le 4|Rm|(\bx,\bt),\end{equation} whenever \begin{equation} (x,t)\in
 M_{\al}, 0<t\le \bt,
 d(x,t)\le d(\bx,\bt)+A|Rm|^{-\frac{1}{2}}(\bx,\bt)\end{equation} }\par {\it Proof
 of Claim 1.}  We construct $(\bx,\bt)$ as a limit of a (finite)
 sequence $(x_k,t_k),$ defined in the following way. Let
 $(x_1,t_1)$ be an arbitrary point, satisfying $0<t_1\le \e^2,
 d(x_1,t_1)<\e, |Rm|(x_1,t_1)\ge \al t^{-1}+\e^{-2}.$ Now if
 $(x_k,t_k)$ is already constructed, and if it can not be taken
 for $(\bx,\bt),$ because there is some $(x,t)$ satisfying (10.2),
 but not (10.1), then take any such $(x,t)$ for
 $(x_{k+1},t_{k+1}).$ Clearly, the sequence, constructed in such a
 way, satisfies $|Rm|(x_k,t_k)\ge 4^{k-1}|Rm|(x_1,t_1)\ge
 4^{k-1}\e^{-2},$ and therefore, $d(x_k,t_k)\le (2A+1)\e.$ Since
 the solution is smooth, the sequence is finite, and its last
 element fits. \par {\bf Claim 2.} {\it For $(\bx,\bt), $ constructed above,
 (10.1) holds whenever \begin{equation} \bt -\frac{1}{2}\al Q^{-1}\le t\le
 \bt, \mbox{dist}_{\bt}(x,\bx)\le
 \frac{1}{10}AQ^{-\frac{1}{2}},\end{equation}
 where $Q=|Rm|(\bx,\bt).$} \par {\it Proof of Claim 2.}  We only need to
 show that if $(x,t)$ satisfies (10.3), then it must satisfy
 (10.1) or (10.2). Since $(\bx,\bt)\in M_{\al},$ we have $Q\ge
 \al\bt^{-1},$ so $\bt -\frac{1}{2}\al Q^{-1}\ge \frac{1}{2}\bt.$
 Hence, if $(x,t)$ does not satisfy (10.1), it definitely belongs
 to $M_{\al}.$ Now by the triangle inequality, $d(x,\bt)\le
 d(\bx,\bt)+\frac{1}{10}AQ^{-\frac{1}{2}}.$ On the other hand, using
 lemma 8.3(b) we see that, as
 $t$ decreases from $\bt$ to $\bt -\frac{1}{2}\al Q^{-1},$ the
 point $x$ can not escape from the ball of radius
 $d(\bx,\bt)+AQ^{-\frac{1}{2}}$ centered at $x_0.$ \par Continuing
 the proof of the theorem, and arguing by contradiction, take
 sequences $\e\to 0,\de\to 0$ and solutions $g_{ij}(t),$ violating
 the statement; by reducing $\e,$ we'll assume that
 \begin{equation}
 |Rm|(x,t)\le  \alpha t^{-1}+2\e^{-2}\  \mbox{whenever}\  0\le t\le \e^2\
 \mbox{and}\
 d(x,t)\le \e \end{equation} Take $A=\frac{1}{100n\e}\to\infty ,$ construct
 $(\bx,\bt),$ and consider solutions
 $u=(4\pi(\bt-t))^{-\frac{n}{2}}e^{-f}$ of the conjugate heat
 equation, starting from $\de$-functions at $(\bx,\bt),$ and
 corresponding nonpositive functions $v.$ \par {\bf Claim 3.}{\it  As $\e
 ,\de\to 0,$ one can find times $\tilde{t}\in[\bt -\frac{1}{2}\al
 Q^{-1},\bt],$ such that the integral $\int_B{v}$ stays bounded
 away from zero, where $B$ is the ball at time $\tilde{t}$ of
 radius $\sqrt{\bt-\tilde{t}}$ centered at $\bx.$ }\par {\it Proof of Claim
 3(sketch).}
 The statement is invariant under scaling, so we can try to take a
  limit of scalings of $g_{ij}(t)$ at points $(\bx,\bt)$ with factors
 $Q.$ If the injectivity radii of the scaled metrics at $(\bx,\bt)$ are bounded
 away from zero, then a smooth limit exists, it is complete and
 has $|Rm|(\bx,\bt)=1$ and $|Rm|(x,t)\le 4$ when $\bt
 -\frac{1}{2}\al \le t\le \bt.$ It is not hard to show that the
 fundamental
 solutions $u$ of the conjugate heat equation converge to such a
 solution on the limit manifold. But on the limit manifold,
 $\int_B{v}$ can not be zero for $\tilde{t}=\bt -\frac{1}{2}\al ,$
 since the evolution equation (9.1) would imply in this case that
 the limit is a gradient shrinking soliton, and this is
 incompatible with $|Rm|(\bx,\bt)=1.$ \par If the injectivity
 radii of the scaled metrics tend to zero, then we can change
 the scaling factor, to make the scaled metrics converge to a flat
 manifold with finite injectivity radius; in this case it is not
 hard to choose $\tilde{t}$ in such a way that $\int_B{v}\to
 -\infty.$
 \par The positive lower bound for $-\int_B{v}$ will be denoted by
 $\be.$ \par Our next goal is to construct an appropriate cut-off
 function. We choose it in the form
 $h(y,t)=\phi(\frac{\tilde{d}(y,t)}{10A\e}),$ where
 $\tilde{d}(y,t)=d(y,t)+200n\sqrt{t},$ and $\phi$ is a smooth
 function of one variable, equal one on $(-\infty,1]$ and
 decreasing to zero on $[1,2].$ Clearly, $h$ vanishes at $t=0$
 outside $B(x_0,20A\e);$ on the other hand, it is equal to one
 near $(\bx,\bt).$ \par Now $\Box h=\frac{1}{10A\e}(d_t-\triangle
 d+\frac{100n}{\sqrt{t}})\phi '-\frac{1}{(10A\e)^2}\phi ''.$ Note
 that $d_t-\triangle t+\frac{100n}{\sqrt{t}}\ge 0$ on the set
 where $\phi '\neq 0 \ \ -$ this follows from the lemma  8.3(a) and
 our assumption (10.4). We may also choose $\phi$ so that $\phi
 ''\ge -10\phi, (\phi ')^2\le 10\phi.$ Now we can compute
 $(\int_M{hu})_t=\int_M{(\Box h)u}\le \frac{1}{(A\e)^2},$ so
 $\int_M{hu}\mid_{t=0}\ge
 \int_M{hu}\mid_{t=\bt}-\frac{\bt}{(A\e)^2}\ge 1-A^{-2}.$ Also, by
 (9.1), $(\int_M{-hv})_t\le \int_M{-(\Box h)v}\le
 \frac{1}{(A\e)^2}\int_M{-hv},$ so by Claim 3, $
 -\int_M{hv}\mid_{t=0}\ge \be \mbox{exp}(-\frac{\bt}{(A\e)^2})\ge
 \be (1-A^{-2}).$ \par From now on we"ll work at $t=0$ only. Let
 $\tilde{u}=hu$ and correspondingly $\tilde{f}=f-\mbox{log}h.$
 Then $$ \be (1-A^{-2})\le -\int_M{hv}=\int_M{[(-2\triangle
 f+|\nabla f|^2-R)\bt -f+n]hu}$$ $$=\int_M{[-\bt |\nabla
 \tilde{f}|^2-\tilde{f}+n]\tilde{u}}+ \int_M{[\bt(|\nabla
 h|^2/h-Rh)-h\mbox{log}h]u}$$
 $$\le\int_M{[-\bt|\nabla\tilde{f}|^2-\tilde{f}-n]\tilde{u}}+A^{-2}+100\e^2$$
 ( Note  that  $\int_M{-uh \log h}$ does   not exceed  the  integral
 of\  $u$\  over \\ $B(x_0,20A\e)\backslash B(x_0,10A\e),$ and
 $\int_{B(x_0,10A\e)}{u}\ge \int_M{\bar{h}u}\ge 1-A^{-2},$\\ where
 $\bar{h}=\phi(\frac{\tilde{d}}{5A\e}))$ \par Now scaling the metric by
 the factor $\frac{1}{2}\bt^{-1}$ and sending $\e,\de$ to zero, we
 get a sequence of metric balls with radii going to infinity, and
 a sequence of compactly supported nonnegative functions
 $u=(2\pi)^{-\frac{n}{2}}e^{-f}$ with $\int{u}\to 1$ and
 $\int{[-\frac{1}{2}|\nabla f|^2-f+n]u}$ bounded away from zero by
 a positive constant. We also have isoperimetric inequalities with
 the constants tending to the euclidean one. This set up is in
 conflict with the Gaussian logarithmic Sobolev inequality, as can
 be seen by using spherical symmetrization. \par {\bf 10.2 Corollary}(from
 the proof) {\it Under the same assumptions, we also have at time $t,
 0<t\le (\e r_0)^2,$ an estimate $Vol B(x,\sqrt{t})\ge c\sqrt{t}^n
 $ for $x\in B(x_0,\e r_0),$ where $c=c(n)$ is a universal
 constant.} \par {\bf 10.3 Theorem.} {\it There exist $\e,\de > 0$ with the following
 property. Suppose $g_{ij}(t)$ is a smooth solution to the Ricci
 flow on $[0,(\e r_0)^2],$ and assume that at $t=0$ we have
 $|Rm|(x)\le r_0^{-2}$ in $B(x_0,r_0),$ and $VolB(x_0,r_0)\ge
 (1-\de)\omega_n r_0^n,$ where $\omega_n$ is the volume of the
 unit ball in $\mathbb{R}^n.$ Then the estimate $|Rm|(x,t)\le (\e
 r_0)^{-2}$ holds whenever $0\le t\le (\e r_0)^2,
 \mbox{dist}_t(x,x_0)<\e r_0.$}\par The proof is a slight
 modification of the proof of theorem 10.1, and is left to the
 reader. A natural question is whether the assumption on the
 volume of the ball is superfluous. \par {\bf 10.4 Corollary}(from 8.2,
 10.1, 10.2) {\it  There exist $\e,\de > 0$ and for any $A>0$ there
 exists $\kappa(A)>0$ with the following property. If $g_{ij}(t)$
 is a smooth solution to the Ricci flow on $[0,(\e r_0)^2],$ such
 that
 at $t=0$ we have $R(x)\ge -r_0^{-2}, Vol(\partial\Omega)^n\ge
 (1-\de)c_nVol(\Omega)^{n-1}$  for any
 $x,\Omega\subset B(x_0,r_0),$ and $(x,t)$ satisfies $A^{-1}(\e
 r_0)^2\le t\le (\e r_0)^2, \mbox{dist}_t(x,x_0)\le Ar_0,$ then
 $g_{ij}(t)$ can not be $\kappa$-collapsed at $(x,t)$ on the
 scales less than $\sqrt{t}.$} \par {\bf 10.5} {\it Remark.} It is
 straightforward to get from 10.1 a version of the Cheeger diffeo
 finiteness theorem for manifolds, satisfying our assumptions on
 scalar curvature and isoperimetric constant on each ball of some
 fixed radius $r_0>0.$ In particular, these assumptions are
 satisfied (for some controllably smaller $r_0$), if we assume a
 lower bound for $\mbox{Ric}$ and an almost euclidean lower bound
 for the volume of the balls of radius $r_0.$ (this follows from
 the
 Levy-Gromov isoperimetric inequality); thus we get one of the
 results of Cheeger and Colding [Ch-Co] under somewhat weaker
 assumptions. \par {\bf 10.6*} Our pseudolocality theorem is similar in
 some respect to the results of Ecker-Huisken [E-Hu] on the mean
 curvature flow.
 \section {Ancient solutions with nonnegative curvature operator and bounded entropy}
 \begin{parag}  In this section we consider smooth solutions to the
 Ricci flow $(g_{ij})_t=-2R_{ij}, -\infty<t\le 0,$ such that for
 each $t$ the metric $g_{ij}(t)$ is  a complete non-flat metric of bounded
 curvature and nonnegative curvature operator. Hamilton discovered
 a remarkable differential Harnack inequality for such solutions;
 we need only its trace version \begin{equation}  R_t+2<X,\nabla
 R>+2\mbox{Ric}(X,X)\ge 0\end{equation} and its corollary, $R_t\ge 0.$ In
 particular, the scalar curvature at some time $t_0\le 0$ controls
 the curvatures for all $t\le t_0.$ \par We impose one more
 requirement on the solutions; namely, we fix some $\kappa >0$ and
 require that $g_{ij}(t)$ be $\kappa$-noncollapsed on all scales
 (the definitions 4.2 and 8.1 are essentially equivalent in this
 case). It is not hard to show that this
 requirement is equivalent to a uniform bound
 on the entropy $S,$ defined as in 5.1 using an arbitrary
 fundamental solution to the conjugate heat
 equation. \end{parag} \begin{parag} Pick an arbitrary point $(p,t_0)$ and define
 $\tilde{V}(\tau), l(q,\tau)$ as in 7.1, for $\tau(t)=t_0-t.$
 Recall that for each $\tau>0$ we can find $q=q(\tau),$ such that
 $l(q,\tau)\le \frac{n}{2}.$ \par {\bf Proposition.}{\it The scalings of
 $g_{ij}(t_0-\tau)$ at $q(\tau)$ with factors $\tau^{-1}$ converge
 along a subsequence of $\tau\to\infty$ to a non-flat gradient shrinking
 soliton. }\par
 {\it Proof (sketch).} It is not hard to deduce from (7.16) that for any
 $\e >0$ one can find $\de >0$ such that both $ l(q,\tau)$ and
 $\tau R(q,t_0-\tau)$  do not exceed $\de ^{-1}$ whenever
 $\frac{1}{2}\bar{\tau}\le \tau\le \bar{\tau}$ and
 $\mbox{dist}_{t_0-\bar{\tau}}^2(q,q(\bar{\tau}))\le
 \e^{-1}\bar{\tau}$ for some $\bar{\tau}>0.$ Therefore, taking
 into account the $\kappa$-noncollapsing assumption, we can take a
 blow-down limit, say $\bar{g}_{ij}(\tau),$ defined for
 $\tau\in(\frac{1}{2},1), (\bar{g}_{ij})_{\tau}=2\bar{R}_{ij}.$
  We may assume also that functions $l$ tend to a locally
  Lipschitz function $\bar{l},$ satisfying (7.13),(7.14) in the
  sense of distributions. Now, since $\tilde{V}(\tau)$ is
  nonincreasing and bounded away from zero (because the scaled
  metrics are not collapsed near $q(\tau)$) the limit function
  $\bar{V}(\tau)$ must be a positive constant; this constant is strictly less
  than $\mbox{lim}_{\tau\to
  0}\tilde{V}(\tau)=(4\pi)^{\frac{n}{2}},$ since $g_{ij}(t)$ is
  not
  flat. Therefore, on the one hand, (7.14) must become an
  equality, hence $\bar{l}$ is smooth, and on the other hand, by
  the description of the equality case in (7.12),
  $\bar{g}_{ij}(\tau)$ must be a gradient shrinking soliton with
  $\bar{R}_{ij}+\bar{\nabla}_i\bar{\nabla}_j
  \bar{l}-\frac{1}{2\tau}\bar{g}_{ij}=0.$ If this soliton is flat,
  then $\bar{l}$ is uniquely determined by the equality in (7.14),
  and it turns out that the value of $\bar{V}$ is exactly
  $(4\pi)^{\frac{n}{2}},$ which was ruled out.\end{parag}
  \begin{cor}
  There is only one  oriented two-dimensional
  solution, satisfying the assumptions stated in 11.1, - the round
  sphere. \end{cor} {\it Proof.}   Hamilton [H 10] proved that round sphere is the
  only
  non-flat oriented nonnegatively curved gradient shrinking soliton
  in dimension two. Thus, the scalings of our ancient solution
  must converge to a round sphere. However, Hamilton [H 10] has also
  shown that an almost round sphere is getting more round under
  Ricci flow, therefore our ancient solution must be round.
  \begin{parag} Recall that for any non-compact complete riemannian
  manifold $M$ of nonnegative Ricci curvature and a point $p\in
  M,$ the function $VolB(p,r)r^{-n}$ is nonincreasing in $r>0;$
  therefore, one can define an asymptotic volume ratio $\cV$ as
  the limit of this function as $r\to\infty.$ \par{\bf
  Proposition.}{\it
  Under assumptions of 11.1, $\cV=0$ for each $t.$} \par {\it Proof.}  Induction on dimension. In
  dimension two the statement is vacuous, as we have just shown.
  Now let $n\ge 3,$  suppose that $\cV>0$ for some $t=t_0,$ and
  consider the asymptotic scalar curvature ratio $\cR=\mbox{lim
  sup}R(x,t_0)d^2(x)$ as $d(x)\to\infty.$ ($d(x)$ denotes
  the distance, at time $t_0,$ from $x$ to some fixed point $x_0$)
  If $\cR=\infty,$ then we can find a sequence of points $x_k$ and
  radii $r_k>0,$ such that $r_k/d(x_k)\to 0,
  R(x_k)r_k^2\to\infty , $ and $R(x)\le 2R(x_k)$ whenever $x\in
  B(x_k,r_k).$  Taking blow-up limit of $g_{ij}(t)$ at $(x_k,t_0)$
  with factors $R(x_k),$ we get a smooth non-flat ancient
  solution, satisfying the assumptions of 11.1, which splits off a
  line (this follows from a standard argument based on the
  Aleksandrov-Toponogov concavity). Thus, we can do dimension
  reduction in this case (cf. [H 4,$\S 22$]). \par If $0<\cR <\infty ,$
  then a similar argument gives a blow-up limit in a ball of
  finite radius; this limit has the structure of a non-flat metric
  cone. This is ruled out by Hamilton's strong maximum principle for
  nonnegative curvature operator. \par
  Finally, if $\cR =0,$ then (in dimensions three and up) it is easy to see that
  the metric is flat. \end{parag} \begin{cor} For every $\e>0$ there
  exists $A< \infty $ with the following property. Suppose we have a
  sequence of ( not necessarily complete) solutions $(g_k)_{ij}(t)$ with
  nonnegative curvature operator, defined on $M_k\times[t_k,0],$
  such that for each $k$ the ball $B(x_k,r_k)$ at time $t=0$ is
  compactly contained in $M_k,$  $\frac{1}{2}R(x,t)\le R(x_k,0)=Q_k$ for all
  $(x,t), t_kQ_k\to -\infty , r_k^2Q_k\to\infty$ as $k\to\infty.$ Then $VolB(x_k,A/\sqrt{Q_k})\le\e(A/\sqrt{Q_k})^n$ at
  $t=0$ if $k$ is large enough.\end{cor} {\it Proof. } Assuming the contrary,
  we may take a blow-up limit (at $(x_k,0)$ with factors $Q_k$)
  and get a non-flat ancient solution with positive asymptotic
  volume ratio at $t=0,$ satisfying the assumptions in 11.1,
  except, may be, the $\kappa$-noncollapsing assumption. But if
  that assumption is violated for each $\kappa>0,$ then $\cV(t)$ is
  not bounded away from zero as $t\to -\infty.$ However, this is
  impossible, because it is easy to see that $\cV(t)$ is
  nonincreasing in $t.$ (Indeed, Ricci flow decreases the volume
  and does not decrease the distances faster than $C\sqrt{R}$ per
  time unit, by  lemma 8.3(b)) Thus,
  $\kappa$-noncollapsing holds for some $\kappa>0,$ and we can
  apply the previous proposition to obtain a contradiction.

  \begin{cor}   For every $w>0$ there exist $B=B(w)<\infty ,
  C=C(w)<\infty , \tau_0=\tau_0(w)>0,$ with the following properties.
  \par (a) Suppose we have a (not necessarily complete) solution
  $g_{ij}(t)$ to the Ricci flow, defined on $M\times [t_0,0],$ so
  that at time $t=0$ the metric ball $B(x_0,r_0)$ is compactly
  contained in $M.$ Suppose that at each time $t, t_0\le t\le 0,$
  the metric $g_{ij}(t)$ has nonnegative curvature operator, and
  $VolB(x_0,r_0)\ge wr_0^n.$ Then we have an estimate $R(x,t)\le
  Cr_0^{-2}+B(t-t_0)^{-1}$ whenever $\mbox{dist}_t(x,x_0)\le
  \frac{1}{4}r_0.$ \par (b) If, rather than assuming a lower bound
  on volume for all $t,$ we
 assume it only for $t=0,$ then the same conclusion holds with
 $-\tau_0r_0^2$ in place of $t_0,$ provided that $-t_0\ge
 \tau_0r_0^2.$\end{cor} {\it Proof.}  By scaling assume $r_0=1.$ (a) Arguing
 by contradiction, consider a sequence of $B,C\to \infty,$ of
 solutions $g_{ij}(t)$ and points $(x,t),$ such that
 $\mbox{dist}_t(x,x_0)\le \frac{1}{4}$ and $ R(x,t)> C+
 B(t-t_0)^{-1}.$ Then, arguing as in the proof of claims 1,2 in
 10.1, we can find a point $(\bx,\bt),$ satisfying
 $\mbox{dist}_{\bt}(\bx,x_0)<\frac{1}{3},
 Q=R(\bx,\bt)>C+B(\bt-t_0)^{-1},$ and such that $R(x',t')\le 2Q$
 whenever $\bt-AQ^{-1}\le t'\le \bt,
 \mbox{dist}_{\bt}(x',\bx)<AQ^{-\frac{1}{2}},$ where $A$ tends to
 infinity with $B,C.$ Applying the previous corollary at
 $(\bx,\bt)$ and using the relative volume comparison, we get a
 contradiction with the assumption involving $w.$ \par (b) Let
 $B(w),C(w)$ be good for (a). We claim that
 $B=B(5^{-n}w),C=C(5^{-n}w)$ are good for (b) , for an appropriate
 $\tau_0(w)>0.$ Indeed, let $g_{ij}(t)$ be a solution with
 nonnegative curvature operator, such that $VolB(x_0,1)\ge w$ at
 $t=0,$ and let $[-\tau ,0]$ be the maximal time interval, where
 the assumption of (a) still holds, with $5^{-n}w$ in place of $w$ and with $-\tau$ in
 place of $t_0.$ Then at time $t=-\tau$ we must have
 $VolB(x_0,1)\le 5^{-n}w.$ On the other hand, from lemma 8.3 (b)
 we see that the ball $B(x_0,\frac{1}{4})$ at time $t=-\tau$
 contains the ball
 $B(x_0,\frac{1}{4}-10(n-1)(\tau\sqrt{C}+2\sqrt{B\tau})) $ at time
 $t=0,$ and the volume of the former is at least as large as the
 volume of the latter. Thus, it is enough to choose
 $\tau_0=\tau_0(w)$ in such a way that the radius of the latter
 ball is $>\frac{1}{5}.$ \par Clearly, the proof also works if
 instead of assuming that curvature operator is nonnegative, we
 assumed that it is bounded below by $-r_0^{-2}$ in the
 (time-dependent) metric ball of radius $r_0,$ centered at $x_0.$
 \begin{parag} From now on we restrict our attention to oriented manifolds of dimension
  three. Under the assumptions in 11.1, the solutions on closed
  manifolds must be quotients of the round $\mathbb{S}^3$ or $\mathbb{S}^2\times\mathbb{R}$ - this is proved
  in the same way as in two dimensions, since the gradient
  shrinking solitons are known from the work of Hamilton [H 1,10]. The
  noncompact solutions are described below.
  \par
  {\bf Theorem.}{\it The set of non-compact ancient solutions , satisfying the
  assumptions of 11.1, is compact modulo scaling. That is , from
  any sequence of such solutions and points $(x_k,0)$ with
  $R(x_k,0)=1,$ we can extract a smoothly converging subsequence,
  and the limit satisfies the same conditions. \par Proof.} To
  ensure a converging subsequence it is enough to show that
  whenever $R(y_k,0)\to\infty,$ the distances at $t=0$ between
  $x_k$ and $y_k$ go to infinity as well. Assume the contrary.
  Define a sequence $z_k$ by the requirement that $z_k$ be the
  closest point to $x_k$ (at $t=0$), satisfying
  $R(z_k,0)\mbox{dist}_0^2(x_k,z_k)=~1.$ We claim that
  $R(z,0)/R(z_k,0)$ is uniformly bounded for
  $z\in B(z_k,2R(z_k,0)^{-\frac{1}{2}}).$ Indeed, otherwise we could
  show, using 11.5 and relative volume comparison in nonnegative
  curvature, that the balls $B(z_k,R(z_k,0)^{-\frac{1}{2}})$ are
  collapsing on the scale of their radii. Therefore, using the
  local derivative estimate, due to W.-X.Shi (see [H 4,$\S 13$]), we get a bound on
  $R_t(z_k,t)$ of the order of $R^2(z_k,0).$ Then we can compare
  $1=R(x_k,0)\ge cR(z_k,-cR^{-1}(z_k,0))\ge cR(z_k,0)$ for some
  small $c>0, $ where the first inequality comes from the Harnack
  inequality, obtained by integrating (11.1). Thus, $R(z_k,0)$ are
  bounded. But now the existence of the sequence $y_k$ at bounded
  distance from $x_k$ implies, via 11.5 and relative volume
  comparison, that balls $B(x_k,c)$ are collapsing - a
  contradiction. \par It remains to show that the limit has
  bounded curvature at $t=0.$ If this was not the case, then we
  could find a sequence $y_i$ going to infinity, such that
  $R(y_i,0)\to\infty$ and $R(y,0)\le 2R(y_i,0)$ for $y\in
  B(y_i,A_iR(y_i,0)^{-\frac{1}{2}}), A_i\to\infty .$ Then the
  limit of scalings at $(y_i,0)$ with factors $R(y_i,0)$ satisfies
  the assumptions in 11.1 and splits off a line. Thus by 11.3 it
  must be a round infinite cylinder. It follows that for large $i$
   each $y_i$ is contained in a round cylindrical "neck" of radius
   $(\frac{1}{2}R(y_i,0))^{-\frac{1}{2}}\to 0,$ - something that can not happen
   in an open manifold of nonnegative curvature.\end{parag} \begin{parag} Fix
   $\e>0.$ Let $g_{ij}(t)$ be an ancient solution on a noncompact
   oriented three-manifold $M,$ satisfying the assumptions in
   11.1. We say that a point $x_0\in M$ is the center of an
   $\e$-neck, if the solution $g_{ij}(t)$ in the set $\{(x,t):
   -(\e Q)^{-1}<t\le 0, \mbox{dist}_0^2(x,x_0)<(\e Q)^{-1}\},$
   where $Q=R(x_0,0),$ is, after scaling with factor $Q,$
   $\e$-close (in some fixed smooth topology) to the corresponding
   subset of the evolving round cylinder, having scalar curvature
   one at $t=0.$ \par {\bf Corollary}  (from theorem 11.7 and its
   proof)
    {\it For any $\e>0$ there exists $C=C(\e,\kappa)>0,$ such that if
    $g_{ij}(t)$ satisfies the assumptions in 11.1, and $M_\e$
    denotes the set of points in $M,$ which are not centers of
    $\e$-necks, then $M_{\e}$ is compact and moreover,
    $\mbox{diam}M_{\e} \le CQ^{-\frac{1}{2}},$ and
    $C^{-1}Q\le R(x,0)\le CQ$ whenever $x\in M_{\e},$ where
    $Q=R(x_0,0)$ for some $x_0\in \partial M_{\e}.$ }\end{parag} {\bf 11.9} {\it Remark.}
     It can be shown that there exists
    $\kappa_0>0,$ such that if an ancient solution on a noncompact
    three-manifold satisfies the assumptions in 11.1 with some
    $\kappa>0,$ then it would satisfy these assumptions with
    $\kappa=\kappa_0.$ This follows from the arguments in 7.3,
    11.2, and the statement (which is not hard to prove) that
    there are no noncompact three-dimensional gradient shrinking
    solitons, satisfying 11.1, other than the round cylinder and
    its $\mathbb{Z}_2$-quotients. \par Furthermore, I believe
    that there is only one (up to scaling) noncompact
    three-dimensional $\kappa$-noncollapsed ancient solution with
    bounded positive curvature - the rotationally symmetric
    gradient steady soliton, studied by R.Bryant. In this
    direction, I have a plausible, but not quite rigorous
    argument, showing that any such ancient solution can be made
    eternal, that is, can be extended for $t\in (-\infty
    ,+\infty);$ also I can prove uniqueness in the class of
    gradient steady solitons. \par {\bf 11.10*} The earlier work on
    ancient solutions  and all that can be found in [H 4, $\S 16-22,25,26$].

    \section {Almost nonnegative curvature in dimension three}
    {\bf 12.1} Let $\phi $ be a decreasing function of one variable, tending to
    zero at infinity. A solution to the Ricci flow is said to have
    $\phi$-almost nonnegative curvature if it satisfies
    $Rm(x,t)\ge -\phi (R(x,t))R(x,t)$ for each $(x,t).$
    \par {\bf Theorem.} {\it Given $\e>0,\kappa>0$ and a function $\phi$ as
    above, one can find $r_0>0$ with the following property. If
    $g_{ij}(t), 0\le t\le T$ is a solution to the Ricci flow on a
    closed three-manifold $M,$ which has $\phi$-almost nonnegative
    curvature and is $\kappa$-noncollapsed on scales $<r_0,$ then
    for any point $(x_0,t_0)$ with $t_0\ge 1$ and $Q=R(x_0,t_0)\ge
    r_0^{-2},$ the solution in
    $\{(x,t):\mbox{dist}^2_{t_0}(x,x_0)<(\e Q)^{-1}, t_0-(\e
    Q)^{-1}\le t\le t_0\}$ is , after scaling by the factor $Q,$
    $\e$-close to the corresponding subset of some ancient
    solution, satisfying the assumptions in 11.1.} \par {\it Proof.}  An
    argument by contradiction. Take a sequence of $r_0$ converging
    to zero, and consider the solutions $g_{ij}(t),$ such that the
    conclusion does not hold for some $(x_0,t_0);$ moreover, by
    tampering with the condition $t_0\ge 1$ a little bit,
    choose among all such $(x_0,t_0),$ in the solution under consideration, the one with nearly the smallest
    curvature $Q.$ (More precisely, we can
    choose $(x_0,t_0)$ in such a way that the conclusion of
    the theorem holds for all $(x,t),$ satisfying $R(x,t)>2Q, t_0-HQ^{-1}\le t\le t_0,$
    where $H\to\infty $ as $r_0\to 0)$  Our goal is to show that the sequence of
    blow-ups of such solutions at such points with factors $Q$
    would converge, along some subsequence of $r_0\to 0,$ to an
    ancient solution, satisfying 11.1.
    \par {\bf Claim 1.} {\it  For each $(\bx,\bt)$ with $t_0-HQ^{-1}\le \bt \le t_0$
    we have $R(x,t)\le 4\bar{Q}$ whenever $\bt-c\bar{Q}^{-1}\le
    t\le \bt$ and $\mbox{dist}_{\bt}(x,\bx)\le
    c\bar{Q}^{-\frac{1}{2}},$ where $\bar{Q}=Q+R(\bx,\bt)$ and
    $c=c(\kappa)>0$ is a small constant.} \par {\it Proof of Claim 1.}  Use the fact
    ( following from the choice of $(x_0,t_0)$ and the description
    of the ancient solutions) that for each $(x,t)$ with
    $R(x,t)>2Q$ and $t_0-HQ^{-1}\le t\le t_0$ we have the estimates $|R_t(x,t)|\le
    CR^2(x,t)$, $
    |\nabla R|(x,t)\le CR^{\frac{3}{2}}(x,t).$ \par {\bf Claim 2.} {\it There
    exists $c=c(\kappa)>0$ and for any $A>0$ there exist
    $D=D(A)<\infty , \rho_0=\rho_0(A)>0,$ with the following
    property. Suppose that $r_0<\rho_0,$ and let $\gamma$ be a
    shortest geodesic with endpoints $\bx,x$ in $g_{ij}(\bt), $
    for some $\bt\in[t_0-HQ^{-1},t_0],$ such that $R(y,\bt)>2Q$ for
    each $y\in\gamma.$ Let $z\in\gamma$ satisfy
    $cR(z,\bt)>R(\bx,\bt)=\bar{Q}.$ Then
    $\mbox{dist}_{\bt}(\bx,z)\ge A\bar{Q}^{-\frac{1}{2}}$ whenever
    $R(x,\bt)\ge D\bar{Q}.$} \par {\it Proof of Claim 2.}  Note that from the choice
    of $(x_0,t_0)$ and the description of the ancient solutions it
    follows that an appropriate parabolic (backward in time)
    neighborhood of a point $y\in\gamma$ at $t=\bt$ is $\e$-close to the
    evolving round cylinder, provided $c^{-1}\bar{Q}\le
    R(y,\bt)\le cR(x,\bt)$ for an appropriate $c=c(\kappa).$ Now
    assume that the conclusion of the claim does not hold, take
    $r_0$ to zero, $R(x,\bt)$ - to infinity, and consider the
    scalings
    around $(\bx,\bt)$ with factors $\bar{Q}.$ We can imagine two
    possibilities for the behavior of the curvature along $\gamma$
    in the scaled metric: either it stays bounded at bounded
    distances from $\bx,$ or not. In the first case we can take a
    limit (for a subsequence) of the scaled metrics along
    $\gamma$ and get a nonnegatively curved almost cylindrical
    metric, with $\gamma$ going to infinity. Clearly, in this case
    the curvature at any point of the limit does not exceed
    $c^{-1};$ therefore, the point $z$ must have escaped to
    infinity, and the conclusion of the claim stands. \par In the
    second case, we can also take a limit along $\gamma;$ it is a
    smooth nonnegatively curved manifold near $\bx$ and has
    cylindrical shape where curvature is large; the radius of the
    cylinder goes to zero as we approach the (first) singular
    point, which is located at finite distance from $\bx;$ the
    region beyond the first singular point will be ignored. Thus,
    at $t=\bt$ we have a metric, which is a  smooth metric of
    nonnegative curvature away from a single singular point $o$.
    Since the metric is cylindrical at points close to $o,$ and
    the radius of the cylinder is at most $\e$ times the distance
    from $o,$ the curvature at $o$ is nonnegative in Aleksandrov
    sense. Thus, the metric near $o$ must be cone-like. In other
    words, the scalings of our metric at points $x_i\to o$ with
    factors $R(x_i,\bt)$ converge to a piece of nonnegatively
    curved non-flat metric cone. Moreover, using claim 1, we see
    that we actually have the convergence of the solutions to the
    Ricci flow on some time interval, and not just metrics at
    $t=\bt.$ Therefore, we get a contradiction with the strong
    maximum principle of Hamilton [H 2].
    \par Now continue the proof of theorem, and recall that we are
    considering scalings at $(x_0,t_0)$ with factor $Q.$ It
    follows from claim 2 that at $t=t_0$ the curvature of the
    scaled metric is bounded at bounded distances from $x_0.$
    This allows us to extract a smooth limit at $t=t_0$ (of
    course, we use the $\kappa$-noncollapsing assumption here). The
    limit has bounded nonnegative curvature (if the curvatures
    were unbounded, we would have a sequence of cylindrical necks
    with radii going to zero in a complete manifold of nonnegative
    curvature). Therefore, by claim 1, we have a limit not only at
    $t=t_0,$ but also in some interval of times smaller than
    $t_0.$ \par We want to show that the limit actually exists for all
    $t<t_0.$ Assume that this is not the case, and let $t'$ be the
    smallest value of time, such that the blow-up limit can be
    taken on $(t',t_0].$ From the differential Harnack inequality
    of Hamilton [H 3] we have an estimate $R_t(x,t)\ge
    -R(x,t)(t-t')^{-1},$ therefore, if $\tilde{Q}$ denotes the
    maximum of scalar curvature at $t=t_0,$ then $R(x,t)\le
    \tilde{Q}\frac{t_0-t'}{t-t'}.$ Hence by lemma 8.3(b)
    $\mbox{dist}_t(x,y)\le \mbox{dist}_{t_0}(x,y)+C$ for all $t,$
    where $C=10n(t_0-t')\sqrt{\tilde{Q}}.$
    \par The next step is needed only if our limit is noncompact.
    In this case there exists $D>0,$ such that for any $y$
    satisfying $d=\mbox{dist}_{t_0}(x_0,y)>D,$ one can find $x$
    satisfying $\mbox{dist}_{t_0}(x,y)=d,
    \mbox{dist}_{t_0}(x,x_0)>\frac{3}{2}d.$ We claim that the
    scalar curvature $R(y,t)$ is uniformly bounded for all such
    $y$ and all $t\in (t',t_0].$ Indeed, if $R(y,t)$ is large,
    then the neighborhood of $(y,t)$ is like in an ancient
    solution; therefore, (long) shortest geodesics $\gamma$ and
    $\gamma_0,$ connecting at time $t$ the point $y$ to $x$ and
    $x_0$ respectively, make the angle close to $0$ or $\pi$ at
    $y;$ the former case is ruled out by the assumptions on
    distances, if $D>10C;$ in the latter case, $x$ and $x_0$ are
    separated at time $t$ by a small neighborhood of $y,$ with
    diameter of order $R(y,t)^{-\frac{1}{2}},$ hence the same must
    be true at time $t_0,$ which is impossible if $R(y,t)$ is too
    large.
    \par Thus we have a uniform bound on curvature  outside a certain compact set, which has uniformly
    bounded diameter for all $t\in (t',t_0].$ Then claim 2 gives a uniform bound on
    curvature everywhere. Hence, by claim 1, we can extend our
    blow-up limit past $t'$ - a contradiction. \par {\bf 12.2
    Theorem.} {\it
    Given a function $\phi$ as above, for any $A>0$ there exists
    $K=K(A)<\infty $ with the following property. Suppose in dimension three we have a solution
    to the Ricci flow with $\phi$-almost nonnegative curvature,
    which satisfies the assumptions of theorem 8.2 with $r_0=1.$
    Then $R(x,1)\le K$ whenever $\mbox{dist}_1(x,x_0)<A.$} \par
    {\it Proof.}  In the first step of the proof we check the following
    \par {\bf Claim.} {\it  There exists $K=K(A)<\infty ,$ such that a point $(x,1)$
    satisfies the conclusion of the previous theorem 12.1 (for
    some fixed small $\e>0$), whenever $R(x,1)>K$ and
    $\mbox{dist}_1(x,x_0)<A.$} \par The proof of this statement
    essentially repeats the proof of the previous theorem (the
    $\kappa$-noncollapsing assumption is ensured by theorem 8.2).
    The only difference is in the beginning. So let us argue by
    contradiction, and suppose we have a sequence of solutions and
    points $x$ with $\mbox{dist}_1(x,x_0)<A$ and
    $R(x,1)\to\infty,$ which do not satisfy the conclusion of
    12.1. Then an argument, similar to the one proving claims 1,2
    in 10.1, delivers points $(\bx,\bt)$ with $\frac{1}{2}\le
    \bt\le 1, \mbox{dist}_{\bt}(\bx,x_0)<2A,$ with
    $Q=R(\bx,\bt)\to\infty ,$ and such that $(x,t)$ satisfies the
    conclusion of 12.1 whenever $R(x,t)>2Q, \bt -DQ^{-1}\le
    t\le\bt, \mbox{dist}_{\bt}(\bx,x)<DQ^{-\frac{1}{2}},$ where
    $D\to\infty.$ (There is a little subtlety here in the
    application of lemma 8.3(b); nevertheless, it works, since we
    need to apply it only when the endpoint other than $x_0$
            either satisfies the conclusion of 12.1, or has scalar
            curvature at most $2Q$)  After such $(\bx,\bt)$ are
            found, the proof of 12.1 applies.
            \par Now, having checked the claim, we can prove the
            theorem by applying the claim 2 of the previous
            theorem to the appropriate segment of the shortest
            geodesic, connecting $x$ and $x_0.$
            \par {\bf 12.3 Theorem.} {\it For any $w>0$ there exist
            $\tau=\tau(w)>0, K=K(w)<\infty, \rho=\rho(w)>0$ with
            the following property. Suppose we have a solution
            $g_{ij}(t)$ to the Ricci flow, defined on $M\times
            [0,T),$ where $M$ is a closed three-manifold, and a point $(x_0,t_0),$
            such that the ball $B(x_0,r_0)$ at $t=t_0$ has volume
            $\ge wr_0^n,$ and sectional curvatures $\ge -r_0^{-2}$
            at each point. Suppose that $g_{ij}(t)$ is
            $\phi$-almost nonnegatively curved for some function
            $\phi$ as above. Then we have an estimate
            $R(x,t)<Kr_0^{-2}$ whenever $t_0\ge 4\tau r_0^2, t\in
            [t_0-\tau r_0^2,t_0], \mbox{dist}_t(x,x_0)\le
            \frac{1}{4}r_0,$ provided that $\phi(r_0^{-2})<\rho.$}
            \par {\it Proof.}  If we knew that sectional curvatures are
            $\ge -r_0^{-2}$ for all $t,$ then we could just apply
            corollary 11.6(b) (with the remark after its proof) and take $\tau(w)=\tau_0(w)/2,
            K(w)=C(w)+2B(w)/\tau_0(w).$ Now fix these values of
            $\tau ,K,$ consider a $\phi$-almost nonnegatively
            curved solution $g_{ij}(t),$  a point $(x_0,t_0)$
            and a radius $r_0>0,$ such that the assumptions of the
            theorem do hold whereas the conclusion does not. We
            may assume that any other point $(x',t')$ and radius
            $r'>0$ with that property has either $t'>t_0$ or
            $t'<t_0 -2\tau r_0^2,$ or $2r'>r_0.$ Our goal is to
            show that $\phi(r_0^{-2})$ is bounded away from zero.
            \par Let $\tau '>0 $ be the largest time interval such
            that $Rm(x,t)\ge -r_0^{-2}$ whenever $t\in[t_0-\tau
            'r_0^2,t_0], \mbox{dist}_t(x,x_0)\le r_0.$ If $\tau
            '\ge 2\tau,$ we are done by corollary 11.6(b). Otherwise,
            by elementary Aleksandrov space theory, we can find at
            time $t'=t_0-\tau 'r_0^2$ a ball $B(x',r')\subset
            B(x_0,r_0)$ with $VolB(x',r')\ge
            \frac{1}{2}\omega_n(r')^n,$ and with radius $r'\ge
            cr_0$ for some small constant $c=c(w)>0.$ By the
            choice of $(x_0,t_0)$ and $r_0,$ the conclusion of our
            theorem holds for $(x',t'),r'.$ Thus we have an
            estimate $R(x,t)\le K(r')^{-2}$ whenever $t\in
            [t'-\tau (r')^2,t'], \mbox{dist}_t(x,x')\le
            \frac{1}{4}r'.$ Now we can apply the previous theorem
            (or rather its scaled version) and get an estimate on
            $R(x,t)$ whenever $t\in [t'-\frac{1}{2}\tau
            (r')^2,t'], \mbox{dist}_t(x',x)\le 10r_0.$ Therefore,
            if $r_0>0$ is small enough, we  have $Rm(x,t)\ge
            -r_0^{-2}$ for those $(x,t),$ which is a contradiction
            to the choice of $\tau '.$
            \par {\bf 12.4 Corollary} (from 12.2 and 12.3) {\it Given a
            function $\phi$ as above, for any $w>0$ one can find
            $\rho>0$ such that if $g_{ij}(t)$ is a $\phi$-almost
            nonnegatively curved solution to the Ricci flow,
            defined on $M\times [0,T),$ where $M$ is a closed three-manifold, and
            if $B(x_0,r_0)$ is a metric ball at time $t_0\ge 1,$
            with $r_0<\rho,$ and such that $\min Rm(x,t_0)$
            over $x\in B(x_0,r_0)$ is equal to $-r_0^{-2},$ then
            $VolB(x_0,r_0)\le wr_0^n.$}
            \section {The global picture of the Ricci flow in
            dimension three}
            \par {\bf 13.1} Let $g_{ij}(t)$ be a smooth solution to the Ricci
            flow on $M\times [1,\infty),$ where $M$ is a closed
            oriented three-manifold. Then, according to [H 6, theorem 4.1], the
            normalized curvatures $\tilde{Rm}(x,t)=tRm(x,t)$ satisfy an estimate of the
            form $\tilde{Rm}(x,t)\ge
            -\phi(\tilde{R}(x,t))\tilde{R}(x,t),$ where $\phi$
            behaves at infinity as $\frac{1}{\mbox{log}}.$ This
            estimate allows us to apply the results 12.3,12.4, and
            obtain the following  \par {\bf Theorem.} {\it For any $w>0$
            there exist $K=K(w)<\infty , \rho=\rho(w)>0,$ such
            that for sufficiently large times $t$ the manifold $M$
            admits a thick-thin decomposition $M=M_{thick}\bigcup
            M_{thin}$ with the following properties. (a) For every
            $x\in M_{thick}$ we have an estimate $|\tilde{Rm}|\le
            K$ in the ball $B(x,\rho(w)\sqrt{t}).$ and the volume
            of this ball is at least $\frac{1}{10}w(\rho(w)\sqrt{t})^n.$  (b) For every
            $y\in M_{thin}$ there exists $r=r(y),
            0<r<\rho(w)\sqrt{t},$ such that for all points in the
            ball $B(y,r)$ we have $Rm\ge -r^{-2},$ and the volume
            of this ball is $<wr^n.$} \par Now the arguments in
            [H 6] show that either $M_{thick}$ is empty for large
            $t,$ or , for an appropriate sequence of $t\to 0$ and
            $w\to 0,$ it converges to a (possibly, disconnected)
            complete hyperbolic manifold of finite volume, whose
            cusps (if there are any) are incompressible in $M.$ On
            the other hand, collapsing with lower curvature bound
            in dimension three is understood well enough to claim
            that, for sufficiently small $w>0,$ $\  M_{thin}$ is homeomorphic to a graph
            manifold. \par The natural questions that remain open
            are whether the normalized curvatures must stay
            bounded as $t\to \infty,$ and whether reducible
            manifolds and manifolds with finite  fundamental group
            can have metrics which evolve smoothly by the Ricci
            flow on the infinite time interval.
            \par {\bf 13.2} Now suppose that $g_{ij}(t)$ is defined on
            $M\times [1,T), T<\infty ,$ and goes singular as $t\to
            T.$ Then using 12.1 we see that, as $t\to T,$ either
            the curvature goes to infinity everywhere, and then
            $M$ is a  quotient of either $\mathbb{S}^3$ or
            $\mathbb{S}^2\times \mathbb{R},$ or the region of
            high curvature in $g_{ij}(t)$ is the union of several
            necks and capped necks, which in the limit turn into
            horns (the horns most likely have finite diameter, but
            at the moment I don't have a proof of that). Then at
            the time $T$ we can replace the tips of the horns by
            smooth caps and continue running the Ricci flow until
            the solution goes singular for the next time, e.t.c.
            It turns out that those tips can be chosen in such a
            way that the need for the surgery will arise only
            finite number of times on every finite time interval.
            The proof of this is in the same spirit,
            as our proof of 12.1; it is technically quite complicated,
            but requires no essentially new ideas. It is likely
            that by passing to the limit in this construction one
            would get a canonically defined Ricci flow through
            singularities, but at the moment I don't have a proof
            of that. (The positive answer to the conjecture in 11.9 on the
            uniqueness of ancient solutions would help here) \par
            Moreover, it can be shown, using an argument based on
            12.2, that every maximal horn at any time $T,$ when the
            solution goes singular, has volume at least $cT^n;$
            this easily implies that the solution is smooth (if nonempty) from
            some finite time on. Thus the topology of the original
            manifold can be reconstructed as a connected sum of
            manifolds, admitting a thick-thin decomposition as in
            13.1, and  quotients of $\mathbb{S}^3$ and
            $\mathbb{S}^2\times\mathbb{R}.$ \par
            {\bf 13.3*} Another  differential-geometric approach to the geometrization conjecture
            is being developed by Anderson [A]; he studies the elliptic
            equations, arising as Euler-Lagrange equations for certain functionals of the
            riemannian metric, perturbing the total scalar curvature
            functional, and one can observe certain parallelism between his work and
            that of Hamilton, especially taking into account that, as we have shown in 1.1,
            Ricci flow is the gradient flow for a functional, that closely resembles
            the total scalar curvature.   \section*{References}

            \ \ \  [A] M.T.Anderson Scalar curvature and geometrization
            conjecture for three-manifolds.  Comparison Geometry
            (Berkeley, 1993-94), MSRI Publ. 30 (1997), 49-82.
            \par
             [B-Em] D.Bakry, M.Emery Diffusions hypercontractives.
             Seminaire de Probabilites XIX, 1983-84, Lecture Notes
             in Math. 1123 (1985), 177-206. \par
              [Cao-C] H.-D. Cao, B.Chow Recent developments on the
              Ricci flow. Bull. AMS 36 (1999), 59-74. \par
               [Ch-Co] J.Cheeger, T.H.Colding On the structure of
               spaces with Ricci curvature bounded below I. Jour.
                Diff. Geom. 46 (1997), 406-480. \par
                [C] B.Chow Entropy estimate for Ricci flow on
                compact two-orbifolds. Jour. Diff. Geom. 33
                (1991), 597-600. \par
               [C-Chu 1] B.Chow, S.-C. Chu A geometric
               interpretation of Hamilton's Harnack inequality for
               the Ricci flow. Math. Res. Let. 2 (1995), 701-718.
               \par [C-Chu 2] B.Chow, S.-C. Chu A geometric
               approach to the linear trace Harnack inequality for
               the Ricci flow. Math. Res. Let. 3 (1996), 549-568.
               \par [D] E.D'Hoker String theory. Quantum fields
               and strings: a course for mathematicians
               (Princeton, 1996-97), 807-1011. \par
               \par [E 1] K.Ecker Logarithmic Sobolev inequalities
               on submanifolds of euclidean space. Jour. Reine
               Angew. Mat. 522 (2000), 105-118. \par
               [E 2] K.Ecker A local monotonicity formula for mean
               curvature flow. Ann. Math. 154 (2001), 503-525.
               \par [E-Hu] K.Ecker, G.Huisken In terior estimates
               for hypersurfaces moving by mean curvature. Invent.
               Math. 105 (1991), 547-569. \par
               [Gaw] K.Gawedzki Lectures on conformal field
               theory. Quantum fields and strings: a course for
               mathematicians (Princeton, 1996-97), 727-805. \par
               [G] L.Gross Logarithmic Sobolev inequalities and
               contractivity properties of semigroups. Dirichlet
               forms (Varenna, 1992) Lecture Notes in Math. 1563
               (1993), 54-88. \par
               [H 1] R.S.Hamilton Three manifolds with positive
               Ricci curvature. Jour. Diff. Geom. 17 (1982),
               255-306. \par
               [H 2] R.S.Hamilton Four manifolds with positive
               curvature operator. Jour. Diff. Geom. 24 (1986),
               153-179. \par
               [H 3] R.S.Hamilton The Harnack estimate for the
               Ricci flow. Jour. Diff. Geom. 37 (1993), 225-243.
               \par [H 4] R.S.Hamilton Formation of singularities
               in the Ricci flow. Surveys in Diff. Geom. 2 (1995),
               7-136. \par
               [H 5] R.S.Hamilton Four-manifolds with positive
               isotropic curvature. Commun. Anal. Geom. 5 (1997),
               1-92. \par
               [H 6] R.S.Hamilton Non-singular solutions of the
               Ricci flow on three-manifolds. Commun. Anal. Geom.
               7 (1999), 695-729. \par
               [H 7] R.S.Hamilton A matrix Harnack estimate for
               the heat equation. Commun. Anal. Geom. 1 (1993),
               113-126. \par
               [H 8] R.S.Hamilton Monotonicity formulas for
               parabolic flows on manifolds. Commun. Anal. Geom. 1
               (1993), 127-137. \par
               [H 9] R.S.Hamilton A compactness property for
               solutions of the Ricci flow. Amer. Jour. Math. 117
               (1995), 545-572. \par
               [H 10] R.S.Hamilton The Ricci flow on surfaces.
               Contemp. Math. 71 (1988), 237-261.
               \par [Hu] G.Huisken Asymptotic behavior for
               singularities of the mean curvature flow. Jour.
               Diff. Geom. 31 (1990), 285-299.
               \par [I] T.Ivey Ricci solitons on compact
               three-manifolds. Diff. Geo. Appl. 3 (1993),
               301-307.
               \par [L-Y] P.Li, S.-T. Yau On the parabolic kernel
               of the Schrodinger operator. Acta Math. 156 (1986),
               153-201. \par
               [Lott] J.Lott Some geometric properties of the
               Bakry-Emery-Ricci tensor. arXiv:math.DG/0211065.

\end{document}